\documentclass[a4paper]{article}

\usepackage{amsmath}
\usepackage{amssymb}
\usepackage[all]{xy}

\newtheorem{theo}{Theorem}[section]
\newtheorem{defi}[theo]{Definition}
\newtheorem{prop}[theo]{Proposition}
\newtheorem{propdef}[theo]{Proposition-Definition}
\newtheorem{coro}[theo]{Corollary}
\newtheorem{lem}[theo]{Lemma}
\newtheorem{rem}[theo]{Remark}
\newtheorem{ex}[theo]{Example}
\newtheorem{conj}[theo]{Conjecture}

\newtheorem{theorem}{Theorem}

\renewcommand{\P}{\mathbf{P}}
\newcommand{\C}{\mathbf{C}}
\newcommand{\R}{\mathbf{R}}
\newcommand{\Q}{\mathbf{Q}}

\renewcommand{\O}{\mathcal{O}}
\renewcommand{\H}{\mathrm{H}}
\newcommand{\Hy}{\mathbf{H}}
\newcommand{\D}{\mathbf{D}}
\newcommand{\h}{\mathrm{h}}

\newcommand{\pr}{\mathrm{pr}}
\newcommand{\sm}{\mathrm{sm}}
\newcommand{\id}{\mathrm{id}}

\newcommand{\ibox}[1]{\mbox{\textit{#1}}}

\newcommand{\halmos}{\hfill $\Box$}
\newcommand{\proof}{\noindent \textbf{Proof.~}}

\renewcommand{\epsilon}{\varepsilon}

\renewcommand{\bar}{\overline}
\renewcommand{\tilde}{\widetilde}
\newcommand{\zero}{^{\circ}}

\newcommand{\sigmat}{\tilde{\Sigma}}

\newcommand{\Crond}{\mathcal{C}}

\author{Thomas Dedieu}
\title{Intrinsic pseudo-volume forms for logarithmic pairs}

\begin{document}

\maketitle

\begin{abstract}
We study an adaptation to the logarithmic case of the
Kobayashi-Eisenman pseudo-volume form, 
or rather an adaptation of its variant defined by Claire Voisin,
for which she replaces holomorphic maps by holomorphic
$K$-correspondences. 
We define an intrinsic logarithmic pseudo-volume form $\Phi_{X,D}$ for
every pair $(X,D)$ consisting of a complex manifold $X$ and a normal
crossing Weil divisor $D$ on $X$, the positive part of which is
reduced. 
We then prove that $\Phi_{X,D}$ is generically non-degenerate
when $X$ is 
projective and $K_X+D$ is ample. This result is analogous to the
classical Kobayashi-Ochiai theorem.
We also show the vanishing of $\Phi_{X,D}$ for a large class
of log-$K$-trivial pairs, which is an important step in the direction
of the Kobayashi conjecture about infinitesimal measure hyperbolicity
in the logarithmic case.  
\end{abstract}

\section*{Introduction}\label{intro}

In the standard non logarithmic case, Kobayashi and Eisenman have
defined an intrinsic 
pseudo-volume form $\Psi_X$ on every complex manifold $X$
(\cite{kobayashi}). The 
definition involves all holomorphic maps from the unit polydisk $\D^n
\subset \C^n$ to $X$. $\Psi_X$ coincides with the Poincar{\'e} hyperbolic
volume form on $X$ when $X$ is a quotient (by a group acting freely and
properly discontinuously) of the unit polydisk $\D^n$. In fact, if $X$
is a smooth curve of genus $g$, then we have the following dichotomy
as a consequence of the Klein-Poincar{\'e} uniformization theorem~: if
$g=0$ or $g=1$, then the universal covering of $X$ is $\P^1$ or $\C$,
and $\Psi_X$ vanishes~; if $g\geqslant 2$, then the universal covering
of $X$ is the unit disk $\D$, and
$\Psi_X$ is induced by the Poincar{\'e} volume form on $\D$. 
For an $n$-dimensional manifold $X$, one expects the situation to
follow the same outline. This is in part 
proved by the Kobayashi-Ochiai theorem (\cite{KO2}), which states that
if $X$ is of general type, then $\Psi_X$ is non degenerate outside a
proper closed algebraic subset of $X$. A variety
$X$ such that $\Psi_X >0$ almost everywhere is said to be
infinitesimal measure
hyperbolic. On the other hand, Kobayashi conjectured that if $X$ is
not of general type, then $\Psi_X$ vanishes on a Zariski open subset
of $X$. The Kobayashi conjecture is proved in the $2$-dimensional
case for algebraic varieties, using the classification of surfaces
(see \cite{GG})~: Green and 
Griffiths show that $\Psi_X = 0$ on a dense Zariski open set when $X$
is covered by abelian varieties, and use the fact that algebraic $K3$
surfaces are swept out by elliptic curves.

It is indeed an important step towards the Kobayashi conjecture to
show that 
if $X$ is a Calabi-Yau variety, then $\Psi_X$ vanishes generically. In 
\cite{Kcorresp}, Claire Voisin defines a new intrinsic pseudo-volume
form $\Phi_{X,an}$, which is a variant of $\Psi_X$,
and for which she is able to show that
a very wide range of Calabi-Yau varieties satisfy the Kobayashi
conjecture
(in fact, she shows that the pseudo-volume form $\Phi_{X,an}$
vanishes on these varieties). 
She also proves a theorem relative to $\Phi_{X,an}$, which
is exactly analogous to the Kobayashi-Ochiai theorem. 
The definition of $\Phi_{X,an}$ is obtained from the definition of $\Psi_X$
by replacing the holomorphic maps from $\D^n$ to $X$ by
$K$-correspondences. A $K$-correspondence between two complex
manifolds $X$ and $Y$ of the same dimension is a closed analytic
subset $\Sigma \subset X \times Y$ satisfying the following
properties~: \\
\emph{(i)} the projections $\Sigma \rightarrow X$ and $\Sigma
\rightarrow Y$ are generically of maximal rank on each irreducible
component of $\Sigma$,\\ 
\emph{(ii)} the first projection $\Sigma \rightarrow X$ is proper,\\ 
\emph{(iii)} for every desingularization $\tau : \sigmat 
\rightarrow \Sigma$, letting $f:= \mathrm{pr}_1 \circ \tau$ and
$g:= \mathrm{pr}_2 \circ \tau$, one has the inequality $R_f \leqslant
R_g$ between the ramification divisors of $f$ and $g$. 
\begin{equation}\label{diag_fg}
\xymatrix{
& \sigmat \ar[d]_{\tau} \ar[ddl]_f \ar[ddr]^g & \\
& \Sigma \ar[dl]^{\pr_1} \ar[dr]_{\pr_2} & \\
X & & Y
}\end{equation}
A $K$-correspondence $\Sigma$ has to be seen as the graph of a
multivalued map between $X$ and $Y$. 
The last condition \emph{(iii)} ensures the existence of a generalized
Jacobian map $(J_{\sigmat})^T : g^* K_Y \rightarrow f^* K_X$. 
This definition of a 
$K$-correspondence, which was introduced in \cite{Kcorresp}, derives
from the notion of $K$-equivalence, for which both projections $\Sigma
\rightarrow X$ and $\Sigma \rightarrow Y$ are birational.

\vspace{0.2cm}
It is nowadays understood that for certain problems, it is more
relevant to consider logarithmic pairs $(X,D)$ rather than simply
considering varieties. In this situation, one replaces the canonical
bundle $K_X$ of $X$ by the log-canonical bundle $K_X(D)$. A first
very classical example of this is given by the study of open
varieties. If $U$ is a complex manifold, such that there exist a compact
variety $X$ and a normal crossing divisor $D \subset X$, such that
$U=X \setminus D$, then the study of the pair $(X,D)$ provides a lot of
information about $U$. 
For example, the Betti cohomology with complex
coefficients of $U$ can be computed as the hypercohomology of
the logarithmic de Rham complex $\Omega_X^\bullet(\log \,D)$, see
\emph{e.g.} \cite{claire}. 
It has also been made clear, that the minimal
model program has to be worked out for pairs, rather than simply for
varieties. But the best clue, showing that it is indeed necessary to
define an intrinsic pseudo-volume form for logarithmic pairs (analogous
to $\Psi_X$), is perhaps the following.

In \cite{campana}, Campana shows that to decompose a
compact K{\"a}hler variety into components of special and hyperbolic
types, one necessarily has to consider fibrations with orbifold
bases. 
By definition, a complex manifold $X$ is of special type if there does
not exist any non trivial meromorphic fibration $X \dashrightarrow Y$
with orbifold base of general type. 
Fano and $K$-trivial manifolds are special, but for every $n>0$ and
$\kappa \in \{-\infty,0,1,\ldots,n-1\}$, there exist $n$-dimensional
manifolds $X$ with $\kappa(X)=\kappa$ that are special. 
If $Y$ is a 
complex manifold, an orbifold structure on $Y$ is the data of a
$\Q$-divisor $\Delta=\sum_j a_j D_j$, where $0 < a_j \leqslant 1$,
$a_j \in \Q$, and the $D_j$ are distinct irreducible divisors on
$Y$. The canonical bundle of the orbifold $(Y/\Delta)$ is the
$\Q$-divisor $K_Y+\Delta$ on $Y$. 
If $X$ and $Y$ are smooth complex varieties, and if $f : X \rightarrow
Y$ is a holomorphic fibration, then the orbifold base of $f$ is
$(Y/\Delta(f))$, where 
$$ \Delta(f) := \sum\nolimits_{D \subset Y} 
\left( 1 - \frac{1}{m(f,D)} \right)\cdot D,$$
$m(f,D)$ being the multiplicity of the fiber of $f$ above the generic
point of $D$.
Campana constructs for every variety $X$ (or rather for every orbifold
$(X/D)$) a functorial fibration $c_X : X \rightarrow C(X)$, the core
of $X$, which is characterized by the fact that the generic fibers are
special, and that the orbifold base is hyperbolic. In addition, he
conjectures 
that the Kobayashi pseudo-metric $d_X$ on $X$ (note that this is not
the same as the Kobayashi-Eisenman pseudo-volume form) is the
pull-back \emph{via} $c_X$ of a pseudo-metric $\delta_X$ on the
orbifold base of the core.

\vspace{0.2cm}
In this paper, we seek the definition of a pseudo-volume form
$\Phi_{X,D}$ on a logarithmic pair $(X,D)$. 
Let $X$ be a complex manifold of dimension $n$, and $D$ be a normal
crossing Weil divisor on $X$, the positive part of which is reduced
(we say that $D$ is a normal crossing divisor if its support has 
normal crossings).
Note that we do not require that $D$ has a non zero positive part.
\begin{theorem}\label{propfondam}
\emph{(i)} 
There exists a logarithmic pseudo-volume form $\Phi_{X,D}$
on the pair $(X,D)$,
\emph{i.e.} a pseudo-metric on the line bundle $\bigwedge^n T_X(-D)$,
satisfying the following functoriality property. 
Let $Y$ be a complex manifold, and 
$\nu : Y \rightarrow X$ be a proper morphism 
with ramification divisor $R$,
such that $\nu^*D-R$ is a normal crossing divisor, 
the positive part of which is reduced.
Then we have 
\begin{equation}\label{functoriality}
\nu^* \Phi_{X,D} = \Phi_{Y,\nu^*D-R}
\end{equation}
(when $\nu$ is not proper, we only get the inequality $\nu^*
\Phi_{X,D} \leqslant \Phi_{Y,\nu^*D-R}$). \\ 
\emph{(ii)}
Let $D$ and $D'$ be two normal crossing Weil divisors on $X$, the
respective positive part of which are reduced.
If $D \leqslant D'$, then $\Phi_{X,D} \leqslant \Phi_{X,D'}$. \\
\emph{(iii)}
If $D=0$, then $\Phi_{X,0}=\Phi_{X,an}$. 
\end{theorem}
This is obtained by following the definition of $\Phi_{X,an}$ in
\cite{Kcorresp}. One 
replaces the holomorphic maps between $\D^n$ and $X$ in the definition
of $\Psi_X$ by log-$K$-correspondences between $(\D^n,\Delta_k)$ and
$(X,D)$, where $\Delta_k$ is the divisor given in $\D^n$ by the equation
$z_{n-k+1} \cdots z_n=0$.
They are closed analytic subsets 
$\Sigma \subset \D^n \times X$, satisfying the following three
properties : 
\emph{(i)} the projections to $X$ and $Y$ are generically of maximal 
rank on each irreducible component of $\Sigma$,
\emph{(ii)} the first projection $\Sigma \rightarrow \D^n$ is proper,
and \emph{(iii)} with the same notations as in (\ref{diag_fg}) above
($\tau : \sigmat \rightarrow \Sigma$ is a desingularization,
$f=\mathrm{pr}_1 \circ \tau$, and $g=\mathrm{pr}_2 \circ \tau$)
\begin{equation}\label{inegalite}
R_f -f^* \Delta_k \leqslant R_g - g^* D.
\end{equation} 
The ramification divisor $R_f$ (resp. $R_g$) is the zero divisor of
the section of $K_{\sigmat} \otimes (f^*K_{\D^n})^{-1}$
(resp. $K_{\sigmat} \otimes (g^*K_X)^{-1}$) given by the Jacobian
map of $f$ (resp. $g$).
Condition \emph{(iii)} ensures the existence
of a generalized Jacobian morphism 
$(J_{\sigmat})^T : g^* K_X(D) \rightarrow f^* K_{\D^n}(\Delta_k)$. 

In \cite{Kcorresp}, Claire Voisin uses local $K$-isocorrespondences to
transport the Poincar{\'e} volume form of $\D^n$ to small open sets on
$X$. Here we use local log-$K$-isocorrespondences to transport the
logarithmic Poincar{\'e} volume form of $(\D^n,\Delta_k)$.
This is done in section \ref{defs}. Log-$K$-correspondences are
defined and studied in section \ref{logKcorresp}.

Section \ref{curvature} is devoted to the proof of the following
result, which generalizes the Kobayashi-Ochiai theorem to our case.
\begin{theorem}\label{KOintro}
Let $(X,D)$ be a pair consisting of a projective $n$-dimensional complex
manifold $X$ and a normal crossing Weil divisor $D$ on $X$,
the effective part of which is reduced and has global normal
crossings. If $K_X+D$ is ample, then $\Phi_{X,D} >0$
(instead of $\Psi_X$, and indeed another motivation for our
construction)
away from a proper closed algebraic subset of
$X$. 
\end{theorem}
As in \cite{Kcorresp}, this is proved by a standard
curvature argument, namely an adaptation to our case of the
Ahlfors-Schwarz lemma, which is more or less an incarnation of the
maximum principle (see \emph{e.g.} \cite{demailly} or
\cite{surveyclaire}). We also use a result of Carlson and Griffiths,
about the existence of metrics with negative Ricci curvature on the
complement of hypersurfaces of projective algebraic manifolds,
enjoying some further properties of K{\"a}hler-Einstein type (see
\cite{Carlson&Griffiths} and \cite{griffiths}).

\vspace{0.2cm}
Let $X$ and $Y$ be two projective complex manifolds of the same
dimension, and assume that $X$ is of general type.
Let $\nu : Y \rightarrow X$ be a dominant morphism.
Using the volume decreasing property $\nu^* \Psi_X \leqslant \Psi_Y$,
it is well known how to obtain an upper bound on $\deg \nu$ from the
classical Kobayashi-Ochiai theorem. However, when $X$ is not of
general type, then we cannot say much about the degree of $\nu$.

One of the major benefits we get by considering $\Phi_{X,D}$ instead of
$\Psi_X$, is that we get much more precise decreasing volume
properties (\emph{cf.} theorem
\ref{propfondam} and section \ref{pseudoVol}).
Indeed, by formula (\ref{functoriality}), where $D=0$,
the standard volume decreasing inequality is even
replaced by an equality for proper morphisms, the ramification
divisors of which have normal crossings.  
Taking the ramification of a morphism $\nu: Y \rightarrow X$ into
consideration gives a greater accuracy in the comparison between the
intrinsic (logarithmic) pseudo-volume forms on $X$ and $Y$
respectively. 

Combining this with theorem \ref{KOintro}, we obtain a way to control
the degree of a dominant morphism $\nu : Y \rightarrow X$ even
if $X$ is not of general type ($X$ and $Y$ are again of the same
dimension). Assume there exists a normal crossing Weil divisor $D$ on
$X$, with positive part reduced and with global normal crossings, 
and such that $K_X+D$ is ample.
One decomposes the ramification divisor $R$ of $\nu$ into
$R = \nu^*D +E$, where $E$ is a Weil divisor of
$Y$, which has a non zero negative part
as soon as $D$ has a non zero positive part. 
Assume that $E$ has normal crossings, 
and that its negative part is reduced. 
Then, by integrating on $Y$ the decreasing volume inequality 
$$ \nu^* \Phi_{X,D} \leqslant \Phi_{Y,-E}, $$
which is an equality if $\nu$ is proper, one gets the inequality
\begin{equation}\label{bound}
\left(\int_X \Phi_{X,D}\right)(\deg \nu) \leqslant \int_Y
\Phi_{Y,-E}. 
\end{equation}
The pseudo-volume forms $\Phi_{X,D}$ and $\Phi_{Y,-E}$ are allowed to
have poles along the positive part of $D$, and the negative part of
$E$ respectively. They are however integrable at the neighbourhood of
these poles, as in the same way the Poincar{\'e} volume form on the punctured
disk $\D \setminus \{0\}$ is integrable at the neighbourhood of $0$.
On the other hand, we have $\int_X \Phi_{X,D}>0$ by theorem
\ref{KOintro}, so inequality (\ref{bound}) yields an upper bound on
$\deg \nu$.

\vspace{0.2cm}
Eventually, we prove our main result in section \ref{log-trivial},
that for many 
log-$K$-trivial pairs, the pseudo-volume form $\Phi_{X,D}$
vanishes. This, of course, has to be seen as a step in the direction
of the logarithmic version of the Kobayashi conjecture. This is stated
as follows.
\begin{theorem}
Let $Y$ be a smooth rationally connected variety, and $(X,D)$ be a pair
such that $X \subset Y$ is a smooth hypersurface, $D$ is reduced and
has normal crossings, and 
$D=X \cap X'$, where $X' \subset Y$ is a reduced hypersurface such that
$$ X + X' \in |-K_Y|. $$
Then the pseudo-volume form $\Phi_{X,D}$ vanishes.
\end{theorem}
Note that by adjunction, the variety $\bar{X}:=X \cup X'$ has trivial
canonical bundle. Using the equality of line bundles $K_X(D) =
\left.(K_{\bar{X}})\right|_X$, we see that the hypotheses impose that
$K_X(D)$ is trivial.
As in \cite{Kcorresp}, the key point in the proof is the production of 
log-$K$-autocorrespondences of the pair $(X,D)$, \emph{i.e.}
correspondences $\Sigma \subset X \times X$ satisfying the equality
$R_f - f^* D = R_g - g^* D$, and with the additional dilating property 
$$f^* \eta_X = \lambda g^* \eta_X,$$
where $|\lambda| \not=1$, and $\eta_X$ is a generator of
$\H^0(X,K_X(D))$. This construction is in some way analogous to the
definition of the multiplication by a complex number on an elliptic
curve, realized as a cubic plane curve.
Another important feature in the proof is the use of pull-backs on
differential forms induced by correspondences, following original
ideas of Mumford (\cite{mumford}, see also chapter 22 in
\cite{claire}).  
One has to be slightly more careful with them than
with pull-backs on cohomology, which are used in \cite{Kcorresp}. 

\vspace{0.2cm}
\textbf{Acknowledgements.} 
I wish to thank Claire Voisin for introducing me to this subject, and
for suggesting this work to me. She helped me to overcome the
traps and difficulties that appeared during its preparation, and
answered all my questions with great clarity.

\section{Log-$K$-correspondences}\label{logKcorresp}

\subsection{Definition and basic properties}\label{BasDefs}

In this section, we define and study the notion of
log-$K$-correspondence. 
This will be used in the next section to define and study properly our
variant of the Kobayashi-Eisenman pseudo-volume form. 
\begin{defi}\label{Kcorresp}
Let $(X,D)$ and $(Y,D')$ be pairs of the same dimension,
\emph{i.e.} $X$ and $Y$ are complex manifolds of dimension $n$, and
$D$ and $D'$ are (not necessarily effective) Weil divisors of $X$ and
$Y$ respectively.  
A log-$K$-correspondence from $(X,D)$ to $(Y,D')$ is a reduced
$n$-dimensional closed analytic subset $\Sigma \subset X \times Y$,
satisfying the three following properties. \\
(i) The projections to $X$ and $Y$ are generically of maximal rank on
each irreducible component of $\Sigma$.\\
(ii) The first projection $\Sigma \rightarrow X$ is proper.\\
(iii) Let $\tau : \sigmat \rightarrow \Sigma$ be a desingularization,
$f= \pr_1 \circ \tau : \sigmat \rightarrow X$, and $g = \pr_2 \circ
\tau : \sigmat \rightarrow Y$. The ramification divisors $R_f$ and
$R_g$ (of $f$ and $g$ respectively) satisfy the inequality
$$ R_f - f^*D \leqslant R_g - g^* D'. $$
\end{defi}
The above notations are summed up in the following commutative
diagram. They will be used very often without further notice in the
end of this text.
\begin{equation}\label{notations}\xymatrix{
& \sigmat \ar[d]_{\tau} \ar[ddl]_f \ar[ddr]^g & \\
& \Sigma \ar[dl]^{\pr_1} \ar[dr]_{\pr_2} & \\
X & & Y
}\end{equation}
Note that if condition \emph{(iii)} is true for one desingularization
of $\Sigma$, then it is true for all desingularizations. Let us now
explain the meaning of this last condition.
The two Jacobian maps $\bigwedge^n df$ and $\bigwedge^n dg$ (or rather
their transpose) give
isomorphisms of line bundles on $\sigmat$
$$
f^* \left( K_X(D) \right) \cong K_{\sigmat}(-R_f+f^*D)
\quad \mbox{and} \quad
g^* \left( K_Y(D') \right) \cong K_{\sigmat}(-R_g+g^*D').
$$
So condition \emph{(iii)} ensures the existence of a holomorphic
(rather than just meromorphic) map
$$
(J_{\sigmat})^T :
g^* \left(K_Y(D')\right) \rightarrow f^* \left(K_X(D)\right),
$$
given by the transpose of the holomorphic map
$$\textstyle{
J_{\sigmat} := \bigwedge^n dg \circ \left( \bigwedge^n df \right)^{-1}
:\ 
f^*\left( \bigwedge^n T_X (-D)\right) 
\rightarrow
g^*\left( \bigwedge^n T_Y (-D') \right) ,
}$$
which we call the generalized logarithmic Jacobian map. 
When $D=0$ and $D'=0$, the notion of log-$K$-correspondence coincides
with the notion of $K$-correspondence introduced in \cite{Kcorresp}.

It is also important to note the following inequality of divisors on
$\sigmat$. We write $D=D_1-D_2$, with both 
$D_1$ and $D_2$ non negative, and similarly $D'=D'_1-D'_2$. If $D_1$
is reduced, then the negative part of $R_f-f^*D$ is the sum of
the reduced divisor $(f^*D_1)_{\mathrm{red}}$ 
and possibly of some $f$-exceptional components contained in
$f^*D_1$. We write this $f$-exceptional sum $E_1$.
In the same way, if $D_1'$ is reduced, the negative part of
$R_g-g^*D'$ writes $(g^*D_1')_{\mathrm{red}} + E_1'$, where $E_1'$ is
a sum of $g$-exceptional components contained in $g^*D_1'$.
So if
both $D_1$ and $D'_1$ are reduced, then condition
\emph{(iii)} implies the inequality
$$ (f^*D_1)_{\mathrm{red}} + E_1 \geqslant 
(g^*D'_1)_{\mathrm{red}} + E_1'. $$
In particular, if $D'$ has a positive part (\emph{i.e.} if $D'_1$ does
not vanish), then $D$ necessarily has a positive part as well.

We shall now describe some enlightening examples.
\begin{ex}\label{antilog}
Let $X$ and $Y$ be complex manifolds of dimension $n$, and $D \subset
X$ be an effective divisor.
A morphism $\varphi : X \rightarrow Y$ such that the ramification
divisor $R_{\varphi}$ contains $D$ (with multiplicities) yields a
morphism of line bundles 
\begin{equation}\label{jaclog}
\varphi^* K_Y \rightarrow K_X(-D). 
\end{equation}
The graph $\Gamma_{\varphi} \subset X \times Y$ is isomorphic to
$X$. It satisfies properties
\emph{(i)} and \emph{(ii)} of definition \ref{Kcorresp}, and with the
notations of (\ref{notations}), one has $R_g - R_f = R_{\varphi}
\geqslant f^* D$. 
So $\Gamma_{\varphi}$ is a (smooth) log-$K$-correspondence between
$(X,-D)$ and $(Y,0)$, and the generalized
logarithmic Jacobian map 
$$ g^* K_Y \rightarrow f^* \left( K_X(-D) \right) $$
identifies to (\ref{jaclog}).
\end{ex}

\begin{ex}\label{log}
Let $(X,D)$ and $(X',D')$ be smooth logarithmic pairs,
where $X$ and $X'$ are complex manifolds of the same
dimension,
and $D \subset X$ and $D' \subset X'$ are effective divisors. One
also usually assumes $D$ and $D'$ to be normal crossing divisors, but
this is not necessary for this example.
A morphism of pairs $\varphi : (X,D) \rightarrow (X',D')$ is a
morphism of complex manifolds $X \rightarrow X'$, such that the
ramification divisor $R_{\varphi}$ contains $\varphi^* D' - D$. In
other words, we require that there exists an effective divisor $R
\subset X$, such that $K_X + D = \varphi^*(K_{X'} + D') + R$ as
divisors on $X$. Such a morphism yields a morphism of line bundles on
$X$
\begin{equation}\label{jacantilog}
\varphi^* \left( K_{X'}(D') \right) \rightarrow K_X(D). 
\end{equation}
Again, the graph $\Gamma_{\varphi} \subset X \times Y$ is isomorphic
to $X$, and 
satisfies both properties \emph{(i)} and \emph{(ii)} of definition
\ref{Kcorresp}, and we have $R_g - R_f = R_{\varphi} \geqslant g^* 
D' - f^* D$. 
So $\Gamma_{\varphi}$ is a log-$K$-correspondence between $(X,D)$ and
$(X',D')$, and the generalized logarithmic Jacobian map
$$ g^* \left( K_{X'}(D') \right) \rightarrow f^* \left( K_X(D)
\right)$$
identifies to (\ref{jacantilog}).
\end{ex}

Finally, the following notion of log-$K$-isocorrespondence will be useful
later.
\begin{defi}\label{logiso}
Let $(X,D)$ and $(Y,D')$ be logarithmic pairs of the same dimension $n$,
and let $\Sigma \subset X \times Y$ be a reduced closed analytic
subset, generically finite over $X$ and $Y$. We let $\tau :
\sigmat \rightarrow \Sigma$ be a desingularization, and use the
notations (\ref{notations}). If both projections $\pr_1$ and $\pr_2$
are proper, and if  
$$ R_f - f^*D = R_g - g^* D', $$
then $\Sigma$ is a log-$K$-isocorrespondence between $(X,D)$ and
$(Y,D')$. 
\end{defi}
Note that under these hypotheses, $\Sigma$ is a log-$K$-correspondence
between $(X,D)$ and $(Y,D')$, and its transpose $\Sigma^T \subset Y
\times X$ is a log-$K$-correspondence between $(Y,D')$ and
$(X,D)$. The generalized logarithmic Jacobian map 
then induces an isomorphism of line bundles on $\sigmat$
$$ g^*\left(K_Y(D')\right) \cong f^*\left(K_X(D)\right).$$

\begin{ex}\label{multpq}
We consider the unit disk $\D$, with maps $f: z\in \D \mapsto z^p \in
\D$ and $g: z \in \D \mapsto z^q \in \D$, where $p$ and $q$ are two
relatively prime integers. Then the diagram
$$\xymatrix{
& \D \ar[dl]_{\displaystyle{z \mapsto z^p}} 
\ar[dr]^{\displaystyle{z \mapsto z^q}} & \\
(\D,\{0\}) && (\D,\{0\})
}$$
yields a log-$K$-autocorrespondence of the pair $(\D,\{0\})$. Indeed,
we have the equality of divisors on $\D$
$$ R_f - f^*\{0\} = (p-1)\{0\} - p\{0\} = -\{0\} = R_g -g^*\{0\}.$$
\end{ex}

\subsection{Composition of log-$K$-correspondences}\label{composition}

We shall now study carefully the composition of two
log-$K$-correspondences. This will allow us in the next section to
prove some properties of volume decreasing type for our logarithmic
pseudo-volume form.

We first need to define the weaker notion of $0$-correspondence, and
to study the composition of two of them.
Let $X$ and $Y$ be two $n$-dimensional complex manifolds. 
\begin{defi}\label{0corresp}
A $0$-correspondence between $X$ and $Y$ is a reduced closed analytic
subset $\Sigma \subset X \times Y$, which is generically finite over
$X$ and $Y$, and such that the first projection $\Sigma \rightarrow X$
is proper. 
\end{defi}
In other words, $\Sigma \subset X \times Y$ has only to satisfy
conditions 
\emph{(i)} and \emph{(ii)} of definition \ref{Kcorresp} to be a
$0$-correspondence. 

Let $Z$ be a third $n$-dimensional complex manifold. We denote by 
$p_{lq}$ the projection of $Z \times X \times Y$ to the $l$-th and
$q$-th factors. 
\begin{propdef}\label{defcomp}
Let $\Sigma \subset X \times Y$ and $\Sigma' \subset Z \times X$ be
two $0$-correspondences. We define 
$\Sigma \circ \Sigma' \subset Z \times Y$ 
as the union of the components of $p_{13}\left( p_{12}^{-1}(\Sigma')
  \cap p_{23}^{-1}(\Sigma) \right)$ on which the projections to $Z$
and $Y$ are generically of maximal rank. Then
$\Sigma \circ \Sigma' \subset Z \times Y$ is a $0$-correspondence. 
\end{propdef}

Before stating the proof of this, let us see on a simple example why
it may be necessary to remove certain irreducible components of
$p_{13}\left( p_{12}^{-1}(\Sigma') \cap p_{23}^{-1}(\Sigma)
\right)$. 
Assume for simplicity that $Z$, $X$ and $Y$ are surfaces. Suppose we
are given two $0$-correspondences $\Sigma \subset X \times Y$ and
$\Sigma' \subset Z \times X$, and that there exist two irreducible
curves $C_Z \subset Z$ and $C_Y \subset Y$, and a point $x_0 \in X$,
such that $\Sigma'$ contains $C_Z \times \{x_0\}$ and $\Sigma$
contains $\{x_0\} \times C_Y$. In other words, $\Sigma'$ contains a
contraction of the curve $C_Z$ to the point $x_0$, and $\Sigma$
contains a blow-up of the point $x_0$ onto the curve $C_Y$. Then
$C_Z \times C_Y$ is an irreducible component of
$p_{13}\left( p_{12}^{-1}(\Sigma') \cap p_{23}^{-1}(\Sigma) \right)$
of dimension 2. It is obvious that this component does not satisfy
condition \emph{(i)} of definition \ref{Kcorresp}. Note that it would
correspond to a blow-up of every point of $C_X$ onto the curve
$C_Z$. 

\vspace{0.1cm} \proof
We have a natural identification between
$p_{12}^{-1}(\Sigma') \cap p_{23}^{-1}(\Sigma)$ and $\Sigma' \times_X
\Sigma$. The first projection $\Sigma' \times_X \Sigma \rightarrow
\Sigma'$ is proper by the stability of properness under base change
($\Sigma \rightarrow X$ is proper by the definition of a 
$0$-correspondence). Since the projection $\Sigma' \rightarrow Z$ is
proper as well, the composition $\Sigma' \times_X \Sigma \rightarrow
\Sigma' \rightarrow Z$ is proper. \emph{A fortiori} $\Sigma \circ
\Sigma'$ satisfies condition \emph{(ii)} of definition
\ref{Kcorresp}. Now condition 
\emph{(i)} of definition \ref{Kcorresp} is clearly satisfied, and one
sees that a component of 
$p_{13}\left( p_{12}^{-1}(\Sigma') \cap p_{23}^{-1}(\Sigma) \right)$
which is generically of maximal rank over both $Z$ and $Y$ is
necessarily of dimension $n$. 

\halmos

We now specify definition \ref{defcomp} to the case of
log-$K$-correspondences.
\begin{prop}\label{logcomp}
Let $D_Z$ and $D_Y$ be Weil divisors of $Z$ and $Y$ respectively, and
$D_X$ and $D'_X$ be Weil divisors of $X$.
Assume $\Sigma'$ is a log-$K$-correspondence between $(Z,D_Z)$ and
$(X,D'_X)$, and $\Sigma$ is a log-$K$-correspondence between $(X,D_X)$
and $(Y,D_Y)$.
If $D'_X \geqslant D_X$, then $\Sigma \circ \Sigma'$ is a
log-$K$-correspondence between $(Z,D_Z)$ and $(Y,D_Y)$.
\end{prop}
We have to prove that the two generalized
logarithmic Jacobian maps 
$J_{\sigmat}^T : g^*\left(K_Y(D_Y)\right) \rightarrow
f^*\left(K_X(D_X)\right)$ and
$J_{\sigmat'}^T : g'^*\left(K_X(D'_X)\right) \rightarrow
f'^*\left(K_Z(D_Z)\right)$
can be composed on $\sigmat' \times_X \sigmat$, to obtain a
generalized logarithmic Jacobian map for 
$\Sigma \circ \Sigma'$. This composition is of course well defined,
since we have 
a morphism $K_X(D_X) \rightarrow K_X(D'_X)$ of line bundles on $X$,
because of the inequality $D_X \leqslant D'_X$.
We just have to show that all this lifts to a desingularization
$\sigmat'' \rightarrow \Sigma \circ \Sigma'$, and that the map we
obtain in this way is
actually the generalized logarithmic Jacobian map of $\Sigma \circ
\Sigma'$. This is a consequence of the 
following lemma \ref{comp}.  

We need some further notations to state the lemma properly. We let
$\tau' : \sigmat' \rightarrow \Sigma'$  
and $\tau : \sigmat \rightarrow \Sigma$ be desingularizations of
$\Sigma'$ and $\Sigma$ 
respectively. We call $f'$ and $g'$
(resp. $f$ and $g$) the maps from $\sigmat'$ (resp. $\sigmat$) to
$Z$ and $X$ (resp. $X$ and $Y$).
Let $\phi$ and $\psi$ be the natural projections from $\sigmat' \times_X
\sigmat$ to $\sigmat'$ and $\sigmat$. We call $\Sigma''$ the
union of the components of $\sigmat' \times_X \sigmat$ on which the
maps $F:=f' \circ \phi$ and $G:= g \circ \psi$ to $Z$ and $Y$ are
generically of maximal rank. We consider a desingularization $\tau'' :
\sigmat'' \rightarrow \Sigma'' \subset \sigmat' \times_X \sigmat$, and
call $\tilde{F}$ and $\tilde{G}$ 
the natural maps from $\sigmat''$ to $Z$ and $Y$.

\begin{lem}\label{comp}
With the above notations (see also the diagram below), we have
$$ R_{\tilde{G}} - R_{\tilde{F}} = (\phi \circ \tau'')^* (R_{g'} -
R_{f'}) + (\psi \circ \tau'')^* (R_{g} - R_{f}) $$
as an equality of divisors on $\sigmat''$.
\end{lem}

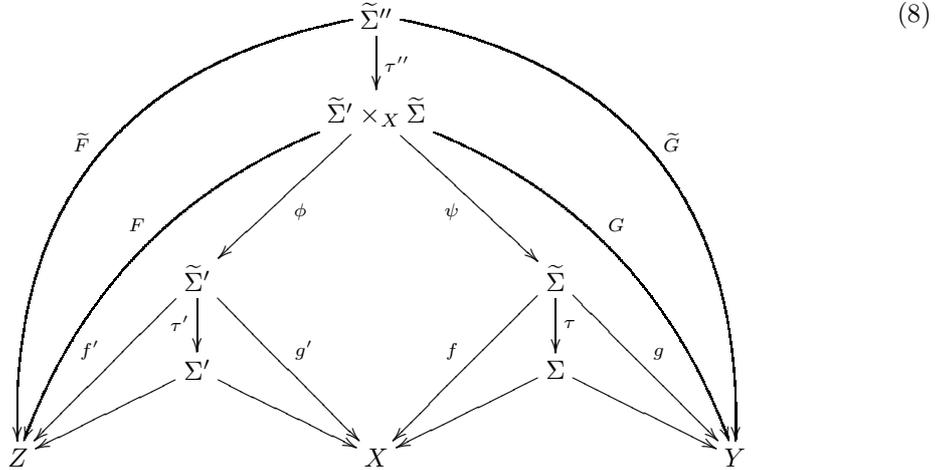
\begin{equation}\label{diagcomp}
\xymatrix@!C@R=20pt{
&& \sigmat'' \ar[1,0]^{\tau''} \ar@/_4pc/[5,-2]_{\tilde{F}}
\ar@/^4pc/[5,2]^{\tilde{G}} && \\
&& \sigmat' \times_X \sigmat \ar[2,-1]^{\phi} \ar[2,1]_{\psi}
\ar@/_2pc/[4,-2]_F \ar@/^2pc/[4,2]^G && \\
\\
& \sigmat' \ar[1,0]_{\tau'} \ar[2,-1]_{f'} \ar[2,1]^{g'} & & \sigmat
\ar[1,0]^{\tau} \ar[2,-1]_{f} \ar[2,1]^{g} \\ 
& \Sigma'\ar[1,-1] \ar[1,1] && \Sigma \ar[1,-1] \ar[1,1] \\
Z && X && Y
}\end{equation}

\proof
Consider $\sigma \in \sigmat''$ and let $z=\tilde{F}(\sigma)$, $x= g'
\circ \phi \circ \tau''(\sigma) = f \circ \psi \circ \tau''(\sigma)$,
and $y= \tilde{G}(\sigma)$.
Let $\omega_z$, $\omega_x$ and $\omega_y$  be holomorphic $n$-forms,
which generate respectively $K_Z$ near $z$, $K_X$ near $x$ and $K_Y$
near $y$. We then have 
$$ {g'}^* \omega_x = \chi' \cdot {f'}^* \omega_z 
\quad \mbox{and} \quad 
g^* \omega_y = \chi \cdot f^* \omega_x,$$
where $\chi'$ is a meromorphic function defined on the inverse image $U'
\subset \sigmat'$ of a neighbourhood of $(z,x) \in Z \times X$, 
with divisor $(\chi') = (R_{g'} - R_{f'}) \cap U'$, and $\chi$ is
similarly a 
meromorphic function defined on $U \subset \sigmat$ with divisor
$(\chi) = (R_{g} - R_{f}) \cap U$. 
Pulling back these equalities on $\sigmat''$ \emph{via} $\phi \circ
\tau''$ and $\psi \circ \tau''$, one gets
$$
{\tau''}^* \phi^* {g'}^* \omega_x = \chi'\circ\phi\circ\tau'' \cdot
{\tau''}^* \phi^* {f'}^* \omega_z  
\quad \mbox{and} \quad 
{\tau''}^* \psi^* g^* \omega_y = \chi\circ\psi\circ\tau'' \cdot
{\tau''}^* \psi^* f^* \omega_x.$$ 
Now since $g' \circ \phi = f \circ \psi$, we have $\phi^* {g'}^* \omega_x
= \psi^* f^* \omega_x$, and therefore
\begin{equation}\label{jac}
\tilde{G}^* \omega_y =  \chi\circ \psi\circ\tau'' \cdot
\chi'\circ\phi\circ\tau'' \cdot \tilde{F}^*\omega_z 
\end{equation}
on the inverse image $U'' \subset \sigmat''$ of a neighbourhood of
$(z,y) \in Z \times Y$. Note that $U''$ contains a neighbourhood of
$\sigma \in \sigmat''$. The meromorphic function $\chi\circ
\psi\circ\tau'' \cdot \chi'\circ\phi\circ\tau''$ has divisor
$$
\left( {\tau''}^* \phi^* (R_{g'} - R_{f'}) + {\tau''}^* \psi^*
  (R_{g} - R_{f}) \right) \cap U''.
$$

\halmos

\vspace{0.2cm}
\noindent \textbf{Proof of proposition \ref{logcomp}.}
Because of proposition \ref{defcomp}, 
to show that $\Sigma \circ \Sigma'$ is a log-$K$-correspondence,
it only remains to prove the inequality
$$ R_{\tilde{F}} - \tilde{F}^* D_Z \leqslant R_{\tilde{G}} -
\tilde{G}^* D_Y. $$
Since $\Sigma'$ and $\Sigma$ are log-$K$-correspondences, we have
$$
R_{g'} - R_{f'} \geqslant g'^*D'_X - f'^*D_Z
\quad \mbox{and} \quad
R_g - R_f \geqslant g^*D_Y - f^* D_X.
$$
By lemma \ref{comp}, this yields
$$ R_{\tilde{G}} - R_{\tilde{F}} \geqslant
{\tau''}^* \phi^* (g'^*D'_X - f'^*D_Z) +
{\tau''}^* \psi^* (g^*D_Y - f^* D_X). $$
On the other hand, since $D'_X \geqslant D_X$, and $g' \circ \phi = f
\circ \psi$, we have
$$ \phi^* g'^*D'_X - \psi^* f^* D_X \geqslant 0, $$
and therefore
$$ R_{\tilde{G}} - R_{\tilde{F}} \geqslant
{\tau''}^* \psi^*g^*D_Y - {\tau''}^* \phi^*f'^*D_Z = 
\tilde{G}^* D_Y - \tilde{F}^* D_Z, $$
which is the desired inequality. 
Eventually, we see from equality (\ref{jac}) that the morphism 
$J_{\sigmat''}^T$ of line bundles on $\sigmat''$ given by the
log-$K$-correspondence $\Sigma \circ \Sigma'$ is obtained as the
composition 
$$\xymatrix@C=5pt{
\tilde{G}^* \left(K_Y(D_Y)\right) \ar[rrr]^{J_{\sigmat''}^T} 
\ar[dr]^(.6){\quad (\tau'' \circ \psi)^* J_{\sigmat}^T} &&& 
\tilde{F}^* \left(K_Z(D_Z)\right) \\
& (\tau'' \circ \psi)^* f^* \left(K_X(D_X)\right) 
\ar@{}[r]|{\subset} & 
(\tau'' \circ \phi)^* g'^* \left(K_X(D'_X)\right) 
\ar[ur]^(.4){(\tau'' \circ \phi)^* J_{\sigmat'}^T\quad }
}$$
of the generalized logarithmic Jacobian maps given by $\Sigma$ and
$\Sigma'$. 

\halmos

\section{Intrinsic logarithmic pseudo-volume forms}\label{defs}

\subsection{The standard Kobayashi-Eisenman pseudo-volume form} 

We first recall the classical definition of the Kobayashi-Eisenman
pseudo-volume form and its fundamental properties.

Let $X$ be an $n$-dimensional complex manifold. The
Kobayashi-Eisenman pseudo-volume form $\Psi_X$ is defined by its
associated Hermitian pseudo-norm on $\bigwedge^n T_X$
$$
\|\xi\|_{\Psi_{X,x}} = \inf \left\{ \lambda>0\ \mbox{\textit{s.t.}}\
  \exists \phi: \D^n \rightarrow 
  X\ \mbox{\textit{with}}\ \phi(0)=x\ \mbox{\textit{and}}\ \lambda \cdot
  J_{\phi}(\partial/\partial 
  t_1 \wedge \cdots \wedge \partial/\partial t_n)=\xi \right\}
$$
for $x \in X$, $\xi \in \bigwedge^n T_{X,x}$, where $\phi$ denotes a
holomorphic map from the unit polydisk in $\C^n$. Note that if $\phi$
is ramified at the origin, then $J_{\phi}(\partial/\partial t_1 \wedge
\cdots \wedge \partial/\partial t_n)=0$, and there does not exist any
$\lambda >0$ such that $\lambda \cdot J_{\phi}(\partial/\partial t_1
\wedge \cdots \wedge \partial/\partial t_n)=\xi$ for $\xi \not= 0$. 

$\Psi_X$ is closely related to the Poincar{\'e} hyperbolic volume form on the
polydisk $\D^n$ 
$$ 
\kappa_n = \frac{i^n}{2^n}
\bigwedge_{1 \leqslant j \leqslant n} \frac{dz_j \wedge
  d\bar{z}_j}{\left( 1-|z_j|^2\right)^2}.
$$
Indeed, $\mathrm{Aut}(\D^n)$ acts transitively on the polydisk and
leaves the 
Poincar{\'e} volume form invariant, and since the latter coincides with
the standard Euclidean volume form at $0$, we find that
$$
\Psi_{X,x} = \inf \left\{ (\phi_b^{-1})^* \kappa_n,\ \phi: \D^n
  \rightarrow X \ \mbox{\textit{s.t.}}\ \phi(b)=x\
  \mbox{\textit{and}}\ \phi\ 
  \mbox{\textit{unramified at}}\ b \right\}, 
$$
where again $\phi$ runs through all holomorphic maps $\D^n \rightarrow
X$ ($\phi_b^{-1}$ is the local inverse of $\phi$ near $b$). One also has
the following result, which is a consequence of
Ahlfors-Schwarz lemma (see section \ref{curvature}).
\begin{theo}[Kobayashi]
If $X$ is isomorphic to the unit polydisk $\D^n \subset \C^n$ (resp. to
the quotient of $\D^n$ by a group acting freely and properly
discontinuously), then the 
Kobayashi pseudo-volume form $\Psi_{X}$ is equal to the Poincar{\'e}
hyperbolic volume form $\kappa_n$ (resp. to the hyperbolic volume form
on the quotient induced by $\kappa_n$).  
\end{theo}

Finally, the following decreasing volume property is a straightforward
consequence of the definition.
If $Y$ is
another smooth manifold of dimension $n$, and if $\phi : X \rightarrow
Y$ is a holomorphic map, then we have the inequality between
pseudo-volume forms on $X$ 
\begin{equation}\label{decVol}
\phi^* \Psi_Y \leqslant \Psi_X.
\end{equation}

There is also a meromorphic version $\tilde{\Psi}_X$ of $\Psi_X$,
introduced by Yau in \cite{yau}. 
It is obtained by considering all meromorphic maps $\phi : \D^n
\dashrightarrow X$ defined near the origin, instead of considering
holomorphic maps as in the definition of $\Psi_X$. It is invariant
under birational maps.

\subsection{The Poincar{\'e} volume form on the punctured
  disk}\label{epointe} 

In this paragraph, we describe the Poincar{\'e} volume form on the
punctured disk $\D \setminus \{0\}$. We will use it later on as a
local model to define intrinsic pseudo-volume forms for logarithmic
pairs. 

The punctured disk $\D \setminus \{0\}$ is the quotient of
$\D$ under the action of the subgroup of $\mathrm{Aut}(\D)$ generated
by a parabolic transformation $g \in \mathrm{Aut}(\D)$. 
In order to compute the Poincar{\'e} volume form of $\D \setminus \{0\}$,
it is however more convenient to see it as a quotient of the Poincar{\'e}
upper half plane $\Hy$. 
\begin{prop}\label{logPoinca}
The Poincar{\'e} volume form of the punctured disk $\D \setminus \{0\}$ is 
$$ \frac{i}{2} \frac{dz \wedge d\bar{z}}{|z|^2 \left( \log |z|^2
  \right)^2}.
$$
\end{prop}

\proof
The punctured disk $\D \setminus \{0\}$ is the quotient of $\Hy$ by
$\left< Z \mapsto Z+1 \right>$. The projection map is
$$ \pi : Z \in \Hy \mapsto z= \exp(2\pi i Z) \in \D \setminus \{0\}.$$
In particular, we have $Z = \log(z)/2\pi i$, and 
$$dZ = \frac{1}{2\pi i} \frac{dz}{z}.$$
Now, we get the Poincar{\'e} volume form on $\D\setminus\{0\}$ by its
expression on $\Hy$
\begin{equation}\label{logPoincaExpr}
\frac{i}{2} \frac{dZ \wedge d\bar{Z}}{|Z-\bar{Z}|^2}
= \frac{i}{2} \frac{dz \wedge d\bar{z}}{|z|^2 \left( \log|z|^2
  \right)^2}.
\end{equation}
\halmos

Proposition \ref{logPoinca} is in fact a particular case of a more general
fact. Let $X$ be a punctured Riemann surface, which 
is universally covered by $\D$. Then every point $x$ in the puncture
corresponds to a subgroup of $\mathrm{Aut}(\D)$ generated by a
parabolic transformation, and there is a neighbourhood of $x$ in $X$
which is isomorphic to the quotient of the circle $\{ \Im(z) >1\}
\subset \Hy$ by $\left<Z \mapsto Z+1\right>$ (see \emph{e.g.}
\cite{FarkasKra}). Therefore, the Poincar{\'e} 
volume form on $X$ is given by (\ref{logPoincaExpr}) around $x$.

\vspace{0.2cm}
The Poincar{\'e} volume form on $\D \setminus \{0\}$ yields a logarithmic
volume form on the pair $(\D,\{0\})$, \emph{i.e.} a volume form with a
pole at $z=0$. 
It is left invariant by the log-$K$-autocorrespondences
\begin{equation}\label{isoDisk}\xymatrix{
& \D \ar[dl]_{\displaystyle{z \mapsto z^p}} 
\ar[dr]^{\displaystyle{z \mapsto z^q}} & \\
(\D,\{0\}) && (\D,\{0\})
}\end{equation}
described in example \ref{multpq}. Of course, it is also left
invariant by the log-$K$-autocorrespondences of $(\D,\{0\})$ given by
the rotations $z \in \D \mapsto e^{i\alpha} \cdot z \in \D$,
$\alpha \in \R$.

\begin{lem}\label{invariante}
Up to multiplication by a constant,
the Poincar{\'e} volume form on $\D \setminus \{0\}$ is
the only logarithmic pseudo-volume form on $(\D,\{0\})$ that is left
invariant by the rotations centred at $0$ and
by the log-$K$-autocor\-res\-pon\-dences (\ref{isoDisk}). 
\end{lem}

\proof
Indeed, let 
$$\frac{i}{2} \alpha(z) \frac{dz\wedge d\bar{z}}{|z|^2}$$
be such a logarithmic volume form on the pair $(\D,\{0\})$.
Since it is invariant under the action of the rotations,
$\alpha(z)$ depends only on $|z|$. On the other hand, letting
$z=z'^p$, we find 
$$\frac{dz}{z} = p \frac{z'^{p-1}}{z'^p} dz'
= p \frac{dz'}{z'}.$$
By invariance under the action of the log-$K$-autocorrespondences
(\ref{isoDisk}), we find that for every $r \in ]0,1[$, and for every
relatively prime integers $p$ and $q$, one has
$$ p^2 \alpha(r^p) = q^2 \alpha(r^q). $$
This implies that there exists $\lambda >0$ such that
$$\frac{i}{2} \alpha(z) \frac{dz\wedge d\bar{z}}{|z|^2}
= \lambda \frac{i}{2} \frac{dz \wedge d\bar{z}}{|z|^2 \left( \log |z|^2
  \right)^2}.
$$
\halmos \pagebreak

This shows in particular that the Poincar{\'e}
logarithmic volume form on $(\D,\{0\})$ is essentially characterized
by its invariance under the action of log-$K$-autocorrespondences of
$(\D,\{0\})$. 
Note that it is locally integrable around $0 \in \D$. 
We will use its $n$-dimensional version given below, 
to define intrinsic logarithmic pseudo-volume forms on general pairs
in the next subsection.

We let $\D^n$ be the unit polydisk, with coordinates
$(z_1,\ldots,z_n)$, and $\Delta_k$ be the divisor given by the
equation $z_{n-k+1}\cdots z_n=0$. Then the pair $(\D^n,\Delta_k)$ is
equipped with the Poincar{\'e} logarithmic volume form
\begin{equation}\label{modele}
\kappa_{n,k} = \left(\frac{i}{2}\right)^n
\left( \bigwedge_{1 \leqslant j \leqslant n-k} 
\frac{dz_j\wedge d\bar{z}_j}{(1-|z_j|^2)^2} \right) \wedge
\left( \bigwedge_{n-k+1 \leqslant j \leqslant n}
\frac{dz_j \wedge d\bar{z}_j}{|z_j|^2 (\log |z_j|^2)^2}\right).
\end{equation}
This is a $C^{\infty}$ logarithmic volume form on $\D^n \setminus
\Delta_k$, and it is singular along $\Delta_k$.

\subsection{Log-$K$-correspondences and intrinsic logarithmic
  pseudo-volume forms}\label{pseudoVol}

In this paragraph, we define the central object of this paper, the
intrinsic pseudo-volume form $\Phi_{X,D}$ of a logarithmic pair
$(X,D)$. 

We first need to introduce the notion of logarithmic pseudo-volume
form. 
\begin{defi}\label{logpform}
Let $(X,D)$ be a pair composed by a complex manifold $X$ of dimension
$n$, and a Weil divisor $D$ of $X$. 
A logarithmic pseudo-volume form on $(X,D)$ is a pseudo-metric on the
line bundle $\bigwedge^n T_X(-D)$.
\end{defi}
Let $\mu$ be a logarithmic pseudo-volume form on $(X,D)$. In case it
is $C^{\infty}$, it writes locally
$$ \mu = \frac{1}{|h|^2} \mu', $$
where $\mu'$ is a $C^{\infty}$ pseudo-volume form, and $h$ is a
meromorphic function with divisor $D$~: if $D=D_1 - D_2$ with $D_1$
and $D_2$ non negative, then $h$ has zeroes exactly along $D_1$, and
poles exactly along $D_2$. It will often be useful to allow $\mu'$ to
have singularities along the positive part $D_1$ of $D$
(\emph{cf.} definition \ref{poincaSing}). 
This is already clear from the expression of the logarithmic Poincar{\'e}
volume form on $(\D,\{0\})$, given in proposition \ref{logPoinca}
above, and more generally from the expression of $\kappa_{n,k}$ in
(\ref{modele}) above.

Now $\Phi_{X,D}$ is defined as follows (we use the notations
introduced in the diagram (\ref{notations})). 
\begin{defi}\label{defphi}
Let $(X,D)$ be a pair composed by an $n$-dimensional complex manifold
$X$ and a normal crossing Weil divisor $D$ of $X$, 
such that the positive part of $D$ is reduced.
For every $x \in X$, we let 
$$\begin{array}{rcl} \Phi_{X,D,x} &=& 
\inf_{0 \leqslant k \leqslant n} \ (\inf \{ (f^* \kappa_{n,k})_{\sigma},\ 
\ibox{where}\ \sigma \in \sigmat,\ \Sigma\
\ibox{log-}K\ibox{-correspondence between} \vspace{0.1cm}\\
&& \hspace{3cm} (\D^n,\Delta_k)\ \ibox{and}\ (X,D),\ 
\ibox{unramified at}\ \sigma,\ \ibox{with}\ g(\sigma)=x
\}).
\end{array}$$
\end{defi}
A log-$K$-correspondence $\Sigma$ between $(\D^n,\Delta_k)$ and
$(X,D)$ is said to be unramified at $\sigma$ if the inequality
of divisors $$ R_f - f^* \Delta_k \leqslant R_g - g^*D $$ is an
equality locally around $\sigma$. 
In this case, the transpose of the generalized
logarithmic Jacobian yields a morphism of line bundles on $\sigmat$
$$ g^* \left( K_X(D) \right) \rightarrow K_{\D^n}(\Delta_k), $$
which is an isomorphism locally around $\sigma$. This authorizes the
identification of $f^* \kappa_{n,k}$ with a Hermitian metric on 
$g^* \left(K_X(D) \right)^{\vee}$ locally around $\sigma$, 
and hence with a logarithmic pseudo-volume element at $x = g(\sigma)$. 

Under the hypotheses of definition \ref{defphi}, there exists an
unramified log-$K$-correponsdence between $(\D^n,\Delta_k)$ 
and $(X,D)$ around every point $x \in X$ 
(the integer $k$ depends on $x$). 
To see this, we write $D = D_1 - D_2$, with $D_1$ and $D_2$ non
negative, and we choose a local holomorphic system of coordinates
$(z_1,\ldots,z_n)$ centred at $x$, and defined on an open set 
$U \subset X$, such
that $D_2$ is given in $U$ by $z_1^{l_1} \cdots z_r^{l_r} = 0$, 
and $D_1$ is given in $U$ by 
$z_{n-k+1} \cdots z_n = 0$ ($r + k \leqslant n$).
Then $D_2$ is the ramification divisor of the morphism 
$U \rightarrow V$ defined by $(z_1,\ldots,z_n) \mapsto 
(z_1^{l_1+1},\ldots,z_r^{l_r+1},z_{r+1},\ldots,z_n)$, where $V$ is an
open subset of $\C^n$. In particular, the graph $\Gamma$ of this
morphism is an unramified log-$K$-correspondence between 
$(X,D_1-D_2)$ and $(V,D_1')$. 
Since $D_1'$ is given in $V$ by the equation
$z_{n-k+1} \cdots z_n = 0$, there exists an unramified
log-$K$-correspondence between $(\D^n,\Delta_k)$ and $(V,D_1')$. 
Then $\Gamma^{T} \circ \Sigma$ is an unramified log-$K$-correspondence 
between $(\D^n,\Delta_k)$ and $(X,D)$ as we wanted. 
Note that $k$ is the number of branches at the point $x$ of the
positive part $D_1$ of $D$.

\begin{prop}\label{ouvertsPoinca}
When $D=0$, one has $\Phi_{X,0}=\Phi_{X,an}$. 
More generally, when $D$ is effective, we have
$$ \left.\Phi_{X,D}\right|_{X\setminus D} = \Phi_{X \setminus D,an}.$$
\end{prop}

\proof
The reason for this is simply that the logarithmic volume form
$\kappa_{n,k}$ on $(\D^n,\Delta_k)$ comes from the Poincar{\'e} volume
form $\kappa_{n,0}$ on $(\D^n,0)$, due to the fact that $\D^n \setminus
\Delta_k$ is a quotient of $\D^n$ by a group acting freely and
properly discontinuously (see paragraph \ref{epointe}). 
Let $\pi_k : \D^n \rightarrow \D^n$ be the projection corresponding to
this quotient. 
One has by definition $\pi_k^* \kappa_{n,k} = \kappa_{n,0}$.

Let $x$ be a point in $X \setminus D$, and 
$\Sigma$ be a log-$K$-correspondence between $(\D^n,\Delta_k)$ and
$(X,D)$, with a point $\sigma \in \sigmat$ above $x$, where $\Sigma$
is unramified. It yields by base change 
$\pi_k : \D^n \rightarrow \D^n$ 
a log-$K$-correspondence $\Sigma'$ between $(\D^n,0)$ and 
$(X\setminus D,0)$. By the same base change, we get a
desingularization $\sigmat'$ of $\Sigma'$.
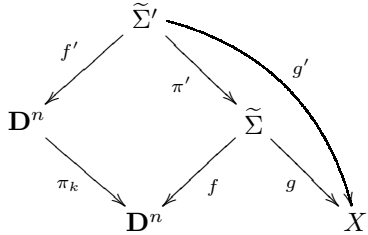
$$\xymatrix{
& \sigmat' \ar[dl]_{f'} \ar[dr]_{\pi'} \ar@/^1.5pc/[2,2]^{g'}  \\
\D^n \ar[dr]_{\pi_k} && \sigmat \ar[dl]^f \ar[dr]_g \\
& \D^n && X
}$$
At every point $\sigma' \in \sigmat'$ above $\sigma$, 
we have $(f'^*\kappa_{n,0})_{\sigma'} = 
(f'^* \pi_k^* \kappa_{n,k})_{\sigma'} =
\pi'^* (f^* \kappa_{n,k})_{\sigma}$. 
In particular, $(f'^*\kappa_{n,0})_{\sigma'}$ and 
$(f^*\kappa_{n,k})_{\sigma}$ yield the same logarithmic volume element
at $x$. 
At $x$, it is therefore sufficient to take the infimum of all $f'^*
\kappa_{n,0}$ in the expression of $\Phi_{X,D}$ in definition
\ref{defphi} (\emph{i.e.} we only consider $k=0$),
with $\Sigma'$ log-$K$-correspondence between
$(\D^n,0)$ and $(X\setminus D,0)$.
By definition, this gives $\Phi_{X \setminus D,an}$ at the point $x$. 

\halmos

\vspace{0.2cm}
When $D=-D_2$ is a negative divisor, 
a similar agument shows
that $\Phi_{X,-D_2}$ can be computed only by looking at
log-$K$-correspondences between $(\D^n,0)$ and $(X,-D_2)$. 
In other words, at a point $x \in X$, $\Phi_{X,-D_2}$ is the infimum
of all $(f^*\kappa_{n,0})_{\sigma}$ where $\sigma \in \sigmat$, and 
$\Sigma$ is a log-$K$-correspondence between $(\D^n,0)$ and
$(X,-D_2)$ as in definition \ref{defphi}. 

\vspace{0.2cm}
One has the following obvious comparison results.
\begin{prop}\label{comparaison}
Let $X$ be a complex manifold. Let $D$ and $D'$ be 
normal crossing Weil divisors of
$X$, such that their respective positive parts are reduced.
If $D \leqslant D'$, then one has 
$$ \Phi_{X,D} \leqslant \Phi_{X,D'}. $$
When $D=0$, one has
$$\Phi_{X,0} = \Phi_{X,an} \leqslant \Psi_X. $$
\end{prop}
The inequality $\Phi_{X,an} \leqslant \Psi_X$
is already contained in \cite{Kcorresp}. It is a simple
consequence of example \ref{antilog}.
The inequality $\Phi_{X,D}\leqslant\Phi_{X,D'}$ when $D\leqslant D'$
comes from the fact that a log-$K$-correspondence between
$(\D^n,\Delta_k)$ and $(X,D')$ is in particular a log-$K$-correspondence
between $(\D^n,\Delta_k)$ and $(X,D)$.

\vspace{0.2cm}
In \cite{Kcorresp}, Claire Voisin shows that if $X$ is a quotient of
$\D^n$ by a group acting freely and properly discontinuously, then
$\Phi_{X,an}$ is simply the Poincar{\'e} volume form on $X$. 
Then by proposition \ref{ouvertsPoinca}, we have 
\begin{equation}\label{logKobayashi}
\Phi_{\D^n,\Delta_k} = \kappa_{n,k}. 
\end{equation}
The proof involves curvature arguments, which are generalized in
section \ref{curvature}. When this is done, we will be able to show
equality (\ref{logKobayashi}) directly (\emph{cf.} theorem
\ref{logKthm}). 

\begin{lem}\label{integrable}
For any pair $(X,D)$ as in definition \ref{defphi}, the logarithmic
pseudo-volume form $\Phi_{X,D}$ is locally integrable.
\end{lem}

\proof
As we already saw, the Poincar{\'e} logarithmic pseudo-volume form 
$$ \frac{i}{2} \frac{dz \wedge d\bar{z}}{|z|^2 (\log|z|^2)^2} $$
on $(\D,\{0\})$ is integrable at the neighbourhood of $0$. As a
consequence, the Poincar{\'e} pseudo-volume form $\kappa_{n,k}$ on 
$(\D^n,\Delta_k)$ is locally integrable as well.

Now for any pair $(X,D)$ as in definition \ref{defphi}, and for any
point $x \in X$, there exists an integer $k$, and a
log-$K$-isocorrespondence $\Sigma$ between $(\D^n,\Delta_k)$ and $(X,D)$,
which is unramified at the neighbourhood of a point $\sigma \in
\sigmat$ over $x$. 
By the expression of $\Phi_{X,D}$ in definition \ref{defphi}, one then
has $(g^* \Phi_{X,D})_{\sigma} \leqslant (f^*
\kappa_{n,k})_{\sigma}$. 
In particular, the growth of $\Phi_{X,D}$ at $x$ is bounded from 
above by the growth of the Poincar{\'e} logarithmic volume form on 
$(\D^n,\Delta_k)$, and hence $\Phi_{X,D}$ is 
locally integrable at $x$.

\halmos

\vspace{0.2cm}
An important feature of these intrinsic pseudo-volume forms is their
decreasing volume properties. In the standard case, this is inequality
(\ref{decVol}). For $\Phi_{X,D}$, we even obtain an equality in the case
of a proper morphism. These properties are obtained by using the study
of the composition of log-$K$-correspondences, which was carried out
in paragraph \ref{composition}.
The main result is the following.

\begin{theo}\label{decvolume}
Let $(X,D)$ and $(Y,D')$ be two pairs composed of a complex
manifold and a normal crossing Weil divisor, 
the positive part of which is reduced.
Assume there exists a
log-$K$-correspondence $\Sigma \subset X \times Y$ between 
$(X,D)$ and $(Y,D')$. We consider a desingularization
$\tau : \sigmat \rightarrow \Sigma$, and use the notations of
(\ref{notations}). 
Then we have the inequality of logarithmic pseudo-volume forms on $\sigmat$ 
$$ g^* \Phi_{Y,D'} \leqslant f^* \Phi_{X, D}. $$
\end{theo}
In the case when $D=D'=0$, this is proved in \cite{Kcorresp}.
As an immediate consequence of theorem \ref{decvolume}, we have the
following. 
\begin{coro}\label{isovolume}
If $\Sigma$ is a log-$K$-isocorrespondence between $(X,D)$ and
$(Y,D')$, 
then we have the equality of logarithmic pseudo-volume forms on
$\sigmat$ 
$$ g^* \Phi_{Y, D'} = f^* \Phi_{X, D}. $$
\end{coro}

\noindent \textbf{Proof.}
Applying theorem \ref{decvolume} to the log-$K$-correspondence
$\Sigma$ on the one hand, and to the log-$K$-correspondence $\Sigma^T$
on the other hand gives both inequalities
$$ g^* \Phi_{Y, D'} \leqslant f^* \Phi_{X, D} 
\quad \mbox{and} \quad
f^* \Phi_{X, D} \leqslant g^* \Phi_{Y, D'}.
$$
One then obviously has the required equality.

\halmos

Applying the former corollary to examples \ref{antilog} and \ref{log},
we get the following.
\begin{coro}\label{pratique}
(i) Let $\pi : X \rightarrow Y$ be a proper morphism, with
ramification divisor $R$ with normal crossings. 
Then we have 
$$\Phi_{X,-R} = \pi^* \Phi_{Y,0}.$$ 
If $\pi$ is not proper, then we only have $\Phi_{X,-R}
\geqslant \pi^* \Phi_{Y,0}$. \\
(ii) Let $\nu : Y \rightarrow X$ be a proper morphism, and $D$ be an
effective, reduced, normal crossing divisor on $X$. We call $D'$ 
the proper transform of $D$ with a reduced scheme
structure. Then $D'$ has normal crossings, and there exist an
effective divisor $R \subset Y$, and an exceptional divisor 
$E \subset Y$, such that $\nu$ ramifies exactly along 
$(\nu^* D -D') -E + R$. If $D'+E-R$ has normal crossings, and if $E$
is reduced, then
$$\nu^* \Phi_{X,D} = \Phi_{Y,D'+E-R}.$$
If $\nu$ is not proper, then we only have $\nu^* 
\Phi_{X,D} \leqslant \Phi_{Y,D'+E-R}$.
\end{coro}

Having the result of proposition \ref{logcomp}, it is fairly easy
to prove the decreasing volume property. The argument is the same as
in the standard case of $\Psi_X$.

\vspace{0.2cm} \noindent \textbf{Proof of theorem \ref{decvolume}.}
Let $\sigma \in \sigmat$, $x=f(\sigma)$, and $y=g(\sigma)$. We have to
show that for every 
$\xi \in f^* \left( \bigwedge^n T_X(-D) \right)_{\sigma}
\subset g^*\left( \bigwedge^n T_Y(-D') \right)_{\sigma}$, one has
the inequality $\Phi_{Y,D'}(g_* \xi) \leqslant \Phi_{X,D}(f_* \xi)$. 

Let $\Sigma' \subset \D^n \times X$ be a log-$K$-correspondence
between $(\D^n, \Delta_k)$ and $(X,D)$, and $\sigma' \in \sigmat'$
such that $g'(\sigma')=x$, and $\Sigma'$ is not ramified at
$\sigma'$. Then there exists 
$\xi' \in g'^*\left( \bigwedge^n T_X(-D) \right)_{\sigma'}$
satisfying $g'_* \xi' = f_* \xi$. By definition, 
$\Phi_{X,D}(f_* \xi)$ is the infimum of the $\kappa_{n,k}(f'_* \xi')$, 
taken on all such $\Sigma'$ and $\xi'$.

Now $\Sigma'' := \Sigma \circ \Sigma' \subset \D^n \times Y$ is a
log-$K$-correspondence between $(\D^n,\Delta_k)$ and $(Y,D')$. It has
a desingularization $\sigmat''$ obtained as in the 
diagram (\ref{diagcomp}). We use the same notations here.
Since $g'(\sigma')=f(\sigma)=x$, there is a point $\sigma'' \in 
\sigmat''$ above both $\sigma$ and $\sigma'$, and there exists
$\xi'' \in (\psi \circ \tau'')^* f^* \left( \bigwedge^n T_X(-D) 
\right)_{\sigma''} = (\phi \circ \tau'')^* g'^*\left( 
\bigwedge^n T_X(-D) \right)_{\sigma''}$, such that
$(\psi \circ \tau'')_*(\xi'') = \xi$ and 
$(\phi \circ \tau'')_*(\xi'') = \xi'$. We thus have 
$\tilde{G}_* \xi'' = g_* \xi$ on the one hand, and 
$\tilde{F}_* \xi'' = f'_* \xi'$, which implies that 
$\kappa_{n,k}(\tilde{F}_* \xi'') = \kappa_{n,k}(f'_* \xi')$
on the other hand. 

Eventually, we have shown the inclusion of sets
\begin{eqnarray*}
&&
\left\{ \kappa_{n,k}(f'_* \xi'),\ \ibox{with}\ \Sigma' \subset 
\D^n \times X,\ \ibox{and}\ g'_* \xi' = f_* \xi \right\} \\
&& \hspace{2cm}
\subset \left\{ \kappa_{n,k}(f''_* \xi''),\ \ibox{with}\ \Sigma'' 
\subset \D^n \times Y,\ \ibox{and}\ g''_* \xi'' = g_* \xi \right\},
\end{eqnarray*}
where $\Sigma'$ (resp. $\Sigma''$) runs through all 
log-$K$-correspondences between $(\D^n,\Delta_k)$ and $(X,D)$ 
(resp. $(Y,D')$). Letting $k$ take any value between $0$ and $n$, and
taking the infima, we get the desired inequality 
$\Phi_{Y,D'}(g_* \xi) \leqslant \Phi_{X,D}(f_* \xi)$.

\halmos

\section{Curvature arguments}\label{curvature}

Let $X$ be a complex manifold.
In the case of the standard Kobayashi-Eisenman pseudo-volume form, one
has the following result, connecting the positivity (or rather the
negativity) of the curvature
of the canonical line bundle of $X$, and the infinitesimal measure
hyperbolicity of $X$ (\emph{i.e.} positivity of $\Psi_X$,
see \cite{demailly}, or \cite{surveyclaire}). 
\begin{theo}[Kobayashi-Ochiai]\label{KO}
If $X$ is a projective complex manifold which is of general type, then
$\Psi_X$ is non-degenerate outside a proper closed algebraic subset of
$X$.
\end{theo}
It is proved in \cite{KO2}, \cite{grif} and for $\tilde{\Psi}_X$
in \cite{yau}. The main 
ingredient in the proof is Ahlfors-Schwarz lemma, which we generalize
in this section.
The converse to theorem \ref{KO} is conjectured to be true in
\cite{kobayashi}. 
\begin{conj}[Kobayashi]\label{kobayashi}
If $X$ is a projective complex manifold which is not of general type,
then $\Psi_X=0$ on a dense Zariski open subset of $X$.
\end{conj}

We study a logarithmic variant of this
conjecture in the logarithmic Calabi-Yau case in section
\ref{log-trivial}. 
In this section, we first establish a generalized Ahlfors-Schwarz
lemma in paragraph \ref{metricNeg}.
Then in paragraph \ref{logGT},
we use it together with a result of Carlson and Griffiths 
to prove a version of the
Kobayashi-Ochiai theorem relative to the logarithmic pseudo-volume
forms $\Phi_{X,D}$. Note that for $D=0$, this is done in
\cite{Kcorresp}. 

\subsection{Metrics with negative curvature on
  $\D^n$}\label{metricNeg} 

We first need to recall a few facts about the logarithmic hyperbolic
volume form $\kappa_{n,k}$ on $(\D^n,\Delta_k)$ and its curvature. 
It is obtained from the Poincar{\'e} metric, which has K{\"a}hler form 
$$
\omega_{n,k} = \frac{i}{2} \left( 
\sum_{1 \leqslant j \leqslant n-k} \frac{dz_j
  \wedge d\bar{z}_j}{\left( 1-|z_j|^2\right)^2}
+ \sum_{n-k+1 \leqslant j \leqslant n} \frac{dz_j \wedge
  d\bar{z}_j}{|z_j|^2 (\log |z_j|^2)^2} \right)
$$
($\kappa_{n,k} = \omega_{n,k}^{n} /n!$, \emph{cf.} (\ref{modele})). It
has Ricci form 
\begin{eqnarray*}
\lefteqn{
-i \partial \bar{\partial} \log \left( \prod_{1 \leqslant j \leqslant n-k}
\frac{1}{(1-z_j \bar{z}_j)^2} 
\prod_{n-k+1 \leqslant j \leqslant n} \frac{1}{z_j \bar{z}_j
  \left(\log(z_j \bar{z}_j)\right)^2 }
\right) }\\
&\hspace{1cm}=& 
2i \left( \sum_{1 \leqslant j \leqslant n-k} \partial \frac{-z_j}{1-z_j
  \bar{z}_j} d\bar{z}_j
+ \sum_{n-k+1 \leqslant j \leqslant n} \partial \frac{1}{\bar{z}_j
  \log(z_j \bar{z}_j)} d\bar{z}_j
\right)
= -4 \omega_{n,k}.
\end{eqnarray*}
This gives the K{\"a}hler-Einstein equation
\begin{equation}\label{KE}
(i \partial \bar{\partial} \log \kappa_{n,k} )^n = 
4^n \omega_{n,k}^n =
4^n n! \kappa_{n,k}.
\end{equation}

Eventually, we need the following definition to state properly the
generalization of Ahlfors-Schwarz lemma in the logarithmic case. 
\begin{defi}\label{poincaSing}
Let $(X,D)$ be a pair composed by a complex manifold $X$, and a Weil
divisor $D=D'-D''$, with $D'$ and $D''$ non negative, and $D'$
reduced and normal crossing.
A logarithmic pseudo-volume form $\mu$ on $(X,D)$ is said to have
singularities of Poincar{\'e} type if it is $C^{\infty}$ on $X\setminus
D'$, and if it is equivalent to a non zero multiple of 
$$ \prod_{n-k+1 \leqslant j \leqslant n} \frac{1}{|z_j|^2
  \left(\log|z_j|^2\right)^2} $$
at the neighbourhood of $D'$, assuming it is given by the equation
$z_{n-k+1} \cdots z_n = 0$.
\end{defi}
Of course, the logarithmic Poincar{\'e} volume form $\kappa_{n,k}$ on
$(\D^n,\Delta_k)$ has singularities of Poincar{\'e} type.

\begin{prop}\label{AS}
Let $(X,D)$ be a pair composed by an $n$-dimensional complex manifold
$X$ and a normal crossing Weil divisor $D$ of $X$, the positive part
of which is reduced.
Let $\Sigma \subset \D^n \times X$ be a log-$K$-correspondence between
$(\D^n,\Delta_k)$ and $(X,D)$, with
desingularization $\tau : \sigmat \rightarrow \Sigma$.
Let $\mu$ be a logarithmic pseudo-volume 
form on $(X,D)$, satisfying the three following properties. \\ 
(a) $i \partial \bar{\partial} \log \mu > 0$. \\
(b) $\mu$ has singularities of Poincar{\'e} type. \\
(c) $(i \partial \bar{\partial} \log \mu)^n \geqslant 4^n n! \mu$.\\
Then one has the inequality of logarithmic pseudo-volume forms on
$\sigmat$
$$ g^* \mu \leqslant f^* \kappa_{n,k}. $$
\end{prop}

Note that if $D$ has a positive part, then $\mu$ has poles, and 
therefore inequality \emph{(c)} cannot be true if $\mu$ has no
singularity, \emph{i.e.} if it writes locally
$\mu'/|h|^2$, where $h$ is a meromorphic function with divisor $D$,
and $\mu'$ is a $C^{\infty}$ pseudo-volume form on $X$. 
In this case, $i \partial \bar{\partial} \log \mu=i \partial
\bar{\partial} \log\mu'$ is a $C^{\infty}$ $(1,1)$-form. In
particular, it has no pole on $X$.

\vspace{0.1cm} 
\noindent \textbf{Proof of proposition \ref{AS}.}
We begin by restricting $\Sigma$ to
$\D^n_{1-\epsilon}:=\D(0,1-\epsilon)^n$ as follows. 
We let $\Sigma_{\epsilon}$ be the inverse image of $\Sigma$ \emph{via}
the map $((1-\epsilon)\id,\id): \D^n \times X \rightarrow \D^n \times
X$. 
This would correspond in the morphism case to the transformation of
$\varphi : \D^n \rightarrow X$ into 
$ \tilde{\varphi}:= \varphi_{|\D^n_{1-\epsilon}} \left(
  (1-\epsilon) \times \cdot\ \right) :
\D^n \rightarrow X. $
Of course one gets from $\sigmat$ a desingularization
$\sigmat_{\epsilon}$ of $\Sigma_{\epsilon}$ in a natural way, with
maps $f_{\epsilon}$ 
and $g_{\epsilon}$ to $\D^n$ and $X$. 

Next, one considers the ratio
$$
\psi_{\epsilon}:= \frac{g_{\epsilon}^* \mu}{f_{\epsilon}^* \kappa_n}.
$$
It is a non-negative $C^{\infty}$-function on
$\sigmat_{\epsilon}$. 
To see why, we first write locally $\mu = \mu'/|h|^2$,
where $h$ is a meromorphic function with divisor $D$
(\emph{i.e.} it has zeroes along the positive part $D_1$ of $D$, and
poles along the negative part $D_2$ of $D$), 
and $\mu'= \chi_{\mu'} (i/2)^n \bigwedge_{1 \leqslant j \leqslant n}
dz_j \wedge d\bar{z}_j$ is a pseudo-volume form on $X$, which is
$C^{\infty}$ on $X \setminus D_1$, and singular along $D_1$, so that
$\mu$ has singularities of Poincar{\'e} type. On the other hand, we let
similarly 
$$\textstyle{ \chi_{n,k} = 
\left(\prod_{1 \leqslant j \leqslant n-k} \left( 1-|t_j|^2\right)^{-2}
  \right)
\left(\prod_{n-k+1 \leqslant j \leqslant n} \left( \log|t_j|^2
  \right)^{-2} \right),}$$ 
so that $\kappa_{n,k} = \chi_{n,k} (i/2)^n \bigwedge_{1 \leqslant j
  \leqslant n} dt_j \wedge d\bar{t}_j/|t_{n-k+1} \cdots t_n|^2$. We
let $h_k(t)=t_{n-k+1} \cdots t_n$ be a holomorphic equation of
$\Delta_k \subset \D^n$.
Then $\psi_{\epsilon}$ writes locally
$$
\psi_{\epsilon}= 
\frac{|s_g|^2 g_{\epsilon}^* \chi_{\mu'}}
{|g_{\epsilon}^*h|^2} \cdot
\frac{|f_{\epsilon}^* h_k|^2}{|s_f|^2 f_{\epsilon}^* \chi_{n,k}},
$$
where $s_f$ and $s_g$ are local analytic equations of the ramification
divisors $R_f$ and $R_g$ respectively.
Now since $\Sigma$ is a log-$K$-correspondence between
$(\D^n,\Delta_k)$ and 
$(X,D)$, one has $R_g -  g^ *D \geqslant R_f - f^* \Delta_k$, and therefore
$s_g \cdot f_{\epsilon}^* h_k/s_f \cdot g_{\epsilon}^*h$ is
$C^{\infty}$ on $X$. 
In addition, since $\mu$ has singularities of Poincar{\'e} type, the ratio 
$g_{\epsilon}^* \chi_{\mu'}/f_{\epsilon}^* \chi_{n,k}$ is $C^{\infty}$
on $X$ as well.

Now, $1/ \chi_{n,k}$ tends to $0$ on the boundary of
$\D^n$, while $ \psi_{\epsilon} \cdot f_{\epsilon}^* \chi_{n,k}$
stays bounded near the boundary of $\sigmat_{\epsilon}$, since 
$\bar{\D}^n_{1-\epsilon}$ is compact and $f$ is
proper.
This gives
$$
\lim_{f_{\epsilon}(x) \rightarrow \partial \D^n} \psi_{\epsilon} (x) =0,
$$ 
and so by properness of $f_{\epsilon}$, the ratio $\psi_{\epsilon}$ has
a maximum on $\sigmat_{\epsilon}$. 
Let $x_0 \in \sigmat_{\epsilon}$ be a point where this maximum is
reached, and write $c:=\psi_{\epsilon}(x_0)$. Then we just have to 
show that $c \leqslant 1$ (we will then get the proposition by letting
$\epsilon$ tend to $0$). 

We now suppose that $c >1$, and show as in the standard proof that
this contradicts hypotheses \emph{(a)}, \emph{(b)} and \emph{(c)}. 
For $ \alpha \in ]1,c[$, define
$$
\sigmat_{\epsilon,\alpha} = \left\{ x \in \sigmat_{\epsilon}\ s.t.\
  \psi_{\epsilon}(x) \geqslant \alpha \right\}.
$$
$\psi_{\epsilon}(x)$ tends to $0$ near the boundary of
$\sigmat_{\epsilon}$, so $\sigmat_{\epsilon,\alpha}$ is compact, and with
smooth boundary for $\alpha$ generic. 
Since 
$i \partial \bar{\partial} \log \chi_{\mu'}=i \partial \bar{\partial}
\log \mu > 0$, and 
$i \partial \bar{\partial} \log \chi_{n,k} = 4 \omega_{n,k} >0$,
$$
\theta :=
g_{\epsilon}^* \left( i \partial \bar{\partial} \log \chi_{\mu'}
\right)^{n-1} 
+ g_{\epsilon}^* \left( i \partial \bar{\partial} \log \chi_{\mu'}
\right)^{n-2}
f_{\epsilon}^* \left( i \partial \bar{\partial} \log \chi_{n,k}
\right)
+ \cdots +
f_{\epsilon}^* \left( i \partial \bar{\partial} \log \chi_{n,k}
\right)^{n-1} 
$$
is a semi-positive $(n-1,n-1)$-form, positive away from the positive
part of $R_f-f^*\Delta_k$.
Note that it is $+\infty$ along the negative part of $R_g-g^*D$,
\emph{i.e.} along the positive part of the reduced divisor
$(g^*D)_{\mathrm{red}}$, where $\mu'$ is singular. 
Now $\psi_{\epsilon}$ has Laplacian
$
i \partial \bar{\partial} \log \psi_{\epsilon} = 
i \partial \bar{\partial} \log (g_{\epsilon}^* \chi_{\mu'})
-i \partial \bar{\partial} \log (f_{\epsilon}^* \chi_{n,k})
$,
and we have 
\begin{eqnarray*}
\left(i \partial \bar{\partial} \log \psi_{\epsilon}\right) \theta 
&=& 
g_{\epsilon}^* \left( i \partial \bar{\partial} \log \chi_{\mu'}
\right)^{n}
- f_{\epsilon}^* \left( i \partial \bar{\partial} \log \chi_{n,k}
\right)^{n} \\
& \geqslant & 
4^n n! \left( g_{\epsilon}^* \mu - f_{\epsilon}^* \kappa_{n,k} \right),
\end{eqnarray*}
where the inequality is given by condition \emph{(c)} on $\mu$, and
by the K{\"a}hler-Einstein equation (\ref{KE}) for the hyperbolic volume
form. 
Now in $\sigmat_{\epsilon,\alpha}$
we have $\psi_{\epsilon}>1$, and therefore
$$
g_{\epsilon}^* \mu - f_{\epsilon}^* \kappa_n \geqslant 0,
$$
with strict inequality away from the positive part of
$R_f-f^*\Delta_k$ (again, this is $+\infty$ along the positive part of 
$(g^*D)_{\mathrm{red}}$). 
Finally the Laplacian $i \partial \bar{\partial} \log \psi_{\epsilon}$
is semi-positive, and positive away from the positive part of
$R_f-f^*\Delta_k$. 

When $x_0$ does not belong to the positive part of
$R_f-f^*\Delta_k$, one can conclude from the maximum principle for 
pluri-subharmonic functions that $\psi_{\epsilon}$ cannot have a
maximum at $x_0$, which is a contradiction, and proves $c \leqslant 1$
as we wanted. Otherwise, one has to apply the following standard
argument. One chooses $m$ satisfying $\log \alpha  < m < \log c$, and then
defines a function 
$$ \mu^+(x) = \max(0,\log \psi_{\epsilon}(x)-m). $$
It is non-negative, positive at $x_0$, and vanishes identically near
the boundary $\partial \sigmat_{\epsilon,\alpha}$. Therefore we have
$$
\int_{ \sigmat_{\epsilon,\alpha}}
\mu^+ \left(i \partial \bar{\partial} \log \psi_{\epsilon}\right)
\theta  > 0.
$$
Note that the form $\left(i \partial \bar{\partial} \log
  \psi_{\epsilon}\right) \theta$ is indeed integrable in
$\sigmat_{\epsilon,\alpha}$, because both $\mu$ and $\kappa_{n,k}$ only
have singularities of Poincar{\'e} type. 
The derivatives of $\mu^+$ are integrable, so we can integrate by
parts the previous integral. This gives
$$
-\int_{ \sigmat_{\epsilon,\alpha}}
\partial \mu^+ \wedge \left(i \bar{\partial} \log \psi_{\epsilon}\right)
\theta  > 0.
$$
Since $\mu^+ = \log \psi_{\epsilon}-m$ when it is non-zero, the former
inequality gives 
$$
\int_{ \sigmat_{\epsilon,\alpha}}
\partial \mu^+ \wedge \left(i \bar{\partial} \log \psi_{\epsilon}\right)
\theta =
\int_{\{ \log \psi_{\epsilon}(x) \geqslant m \}}
i \left(\partial \log \psi_{\epsilon}\right) \wedge\left(
  \bar{\partial} \log \psi_{\epsilon}\right) \theta  < 0,
$$
which is a contradiction, since the right-hand side integral is
obviously positive. 

\halmos

\vspace{0.2cm}
As a first application of this result, we can show directly that
$\Phi_{\D^n,\Delta_k}$ is indeed given by the hyperbolic logarithmic
Poincar{\'e} volume form (see proposition \ref{ouvertsPoinca} as well).
\begin{theo}\label{logKthm}
For $0 \leqslant k \leqslant n$, we have the equality of logarithmic
volume forms on $(\D^n,\Delta_k)$ 
$$ \Phi_{\D^n,\Delta_k} = \kappa_{n,k}. $$
\end{theo}

\proof
The diagonal in $\D^n \times \D^n$ is a log-$K$-correspondence between
$(\D^n,\Delta_k)$ and itself. By definition \ref{defphi}, we thus have
$\Phi_{\D^n,\Delta_k} \leqslant \kappa_{n,k}$.

On the other hand, let $\Sigma$ be a log-$K$-correspondence between
$(\D^n,\Delta_p)$ and $(\D^n,\Delta_k)$, with $0 \leqslant p \leqslant
n$. Then by proposition \ref{AS}, and using the standard notations of
(\ref{notations}), we have $g^* \kappa_{n,k} \leqslant f^*
\kappa_{n,p}$. 
This implies that $\kappa_{n,k} \leqslant \Phi_{\D^n,\Delta_k}$.

\halmos

\subsection{Mappings onto pairs with positive logarithmic canonical
  bundle}\label{logGT} 

Before we state our generalization of the Kobayashi-Ochiai theorem
\ref{KO}, 
let us first recall that a variety $X$ is said to be 
of general type if the canonical bundle $K_X$ is big.
A line bundle $L$ on $X$ is said to be big if
it has maximal Iitaka-Kodaira dimension $\kappa(X,L)=\dim X$. This is
equivalent to the fact, that the
image of the rational map associated to the linear system $|mL|$ is
of maximal dimension $\dim X$, for $m$ big enough and divisible
enough (see \emph{e.g.} \cite{ueno}). 
Given a Weil divisor $D$ on $X$, the condition corresponding to
the fact that $X$ is of general type is in our case the bigness of
$K_X + D$. 
However, we shall need some slightly stronger hypotheses to prove our
result. 
\begin{theo}\label{logKO}
Let $(X,D)$ be a pair composed by a projective $n$-dimensional complex
manifold $X$, and a normal crossing Weil divisor $D$ of $X$, the
positive part of which is reduced and has global normal crossings. 
If $K_X+D$ is ample, then $\Phi_{X,D} >0$ away from a proper closed
algebraic subset of $X$. 
\end{theo}

As an important consequence of this result, we can bound from above
the degree of a morphism of logarithmic pairs, which is onto a pair
with positive logarithmic canonical bundle.

\begin{coro}\label{logKOgenuine}
Let $(X,D)$ be a pair composed by a projective manifold $X$, and a
normal crossing Weil divisor $D$, the positive part of which is 
reduced, and has global normal crossings. 
We assume that $K_X+D$ is ample. 
Let $(Y,D')$ be another logarithmic pair. We assume that
$D'$ has normal crossings, and that its positive part is reduced.
For every dominant morphism
$\phi : (Y,D') \rightarrow (X,D)$, we have
$$ \deg \phi \leqslant \frac{\textstyle{\int_Y
    \Phi_{Y,D'}}}{\textstyle{\int_X \Phi_{X,D}}}. $$
\end{coro}

\proof
The morphism $\phi$ induces a log-$K$-correspondence between $(Y,D')$
and $(X,D)$ (\emph{cf.} example \ref{log}). 
By decreasing volume property, we thus have 
\begin{equation}\label{inegLog}
 \phi^* \Phi_{X,D} \leqslant \Phi_{Y,D'}. 
\end{equation}
Since for $0 \leqslant k \leqslant n$, the volume form $\kappa_{n,k}$
is integrable at the neighbourhood of $\Delta_{k}$, and since $X$ and
$Y$ are compact, definition
\ref{defphi} implies that the pseudo-volume forms $\Phi_{X,D}$ and
$\Phi_{Y,D'}$ are integrable on $X$ and $Y$ respectively.
In particular, (\ref{inegLog}) yields 
$$ \left(\deg \phi\right) \textstyle{\int_X \Phi_{X,D}}
\leqslant \textstyle{\int_Y \Phi_{Y,D'}}. $$
Now, by theorem \ref{logKO}, the integral $\textstyle{\int_X
\Phi_{X,D}}$ is positive, so $\deg \phi$ is bounded from above.

\halmos

\vspace{0.2cm}
Along with Ahlfors-Schwarz lemma \ref{AS}, the main ingredient in the
proof of theorem \ref{logKO} above is the following result.

\begin{lem}[Carlson-Griffiths, \cite{Carlson&Griffiths}]\label{CG}
Let $(X,D)$ be a pair as in theorem \ref{logKO}. Then there exists a
logarithmic volume form with Poincar{\'e} singularities $\mu$ on
$(X,D)$, such that $i \partial \bar{\partial} \log \mu >0$ and $(i
\partial \bar{\partial} \log \mu)^n \geqslant \mu$.
\end{lem}

For a complete proof of this, we refer to \cite{griffiths},
proposition 2.17. Let us still explain how this pseudo-volume form is
constructed. We choose a Hermitian metric $h_0$ on $K_X$. Equivalently,
$h_0^{-1}$ is a $C^{\infty}$ volume form $\mu_X$ on $X$. We assume
for simplicity that $D$ is effective, the general case follows
easily. Then $D$ is assumed to have global normal crossings, so it
writes $D= D_1 + \cdots + D_p$, where the $D_j$ are smooth divisors on
$X$, meeting transversaly. We choose a Hermitian metric $h_j$ for
every line bundle $\O_X(D_j)$, and we let $s_j$ be a global section of
this line bundle, with zero divisor $D_j$. Then for $\alpha >0$
sufficiently small, the pseudo-volume form
$$ \mu_{\alpha} :=
\frac{\mu_X}{\prod_{j=1}^p h_j(s_j) \left(\log(\alpha
    h_j(s_j))\right)^2 }
$$
satisfies the required properties.

Before we continue, we want to underline that this result is used
as a first step in the construction of complete K{\"a}hler-Einstein
metrics with negative Ricci curvature on the complement of
hypersurfaces of projective algebraic manifolds (see
\cite{kobayashiKE} and \cite{TianYau}).

\vspace{0.2cm}\noindent \textbf{Proof of theorem \ref{logKO}.}
By lemma \ref{CG}, there exists a logarithmic pseudo-volume form $\mu$
on $(X,D)$, with Poincar{\'e} singularities, and such that $i \partial
\bar{\partial} \log \mu >0$, and $(i \partial \bar{\partial} \log
\mu)^n \geqslant \mu$. Up to a  rescaling, one can assume that $\mu$
satisfies hypothesis \emph{(c)} of lemma \ref{AS}.
Then if $0 \leqslant k \leqslant n$, and for every
log-$K$-correspondence $\Sigma$ between $(\D^n,\Delta_k)$ and $(X,D)$,
we have the inequality of pseudo-volume forms on a desingularization
$\sigmat$
$$ g^* \mu \leqslant f^* \kappa_{n,k}.$$
It thus follows from the definition \ref{defphi} of $\Phi_{X,D}$ that
$\Phi_{X,D} \geqslant g^* \mu$, which implies that $\Phi_{X,D}$ is
positive on a Zariski dense open subset of $X$, since $\mu$ is a true
logarithmic volume form on $(X,D)$.

\halmos

\section{Log-$K$-autocorrespondences on log-$K$-trivial pairs}\label{log-trivial}

In this section, we prove that for many pairs $(X,D)$, where $X$ is a
complex manifold, and $D$ an effective divisor on $X$, which is
reduced and has normal crossings, such that the
line bundle $K_X(D)$ is trivial, the pseudo-volume form $\Phi_{X,D}$
vanishes. This can be interpreted as a special case of the Kobayashi
conjecture \ref{kobayashi} in the logarithmic case. Again, the case
$D=0$ is handled in \cite{Kcorresp}. However, one has to be slightly
more careful in our case for the proof. This is due to the fact that
we need to pull-back differential forms \emph{via} correspondences,
instead of cohomology classes.

\subsection{Log-$K$-autocorrespondences and the Kobayashi
  conjecture}\label{struct_preuve}

In view of theorem \ref{logKO}, which has to be seen as a logarithmic
version of the Kobayashi-Ochiai theorem for $\Phi_{X,D}$, the Kobayashi
conjecture generalizes as follows in the logarithmic case.
\begin{conj}
Let $(X,D)$ be a pair composed by a projective manifold $X$ and a
normal crossing Weil divisor $D \subset X$, the positive part of 
which is reduced. If $(X,D)$ is not of log-general type
(\emph{i.e.} if $K_X+D$ is not big), then $\Phi_{X,D}$ vanishes on a
dense Zariski open subset of $X$.
\end{conj}

This section is devoted to the proof of the following result, which
goes in the direction of this conjecture.
\begin{theo}\label{annulation}
Let $(X,D)$ be a pair consisting of a smooth projective variety $X$, and
an effective divisor $D \subset X$, which is reduced and has normal
crossings, such that $K_X(D)$ is trivial. 
Assume that there exists a smooth, rationally connected variety $Y$, such
that $X$ can be realized as a hypersurface $X \subset Y$, $D=X \cap
X'$, where $X' \subset Y$ is a reduced 
hypersurface such that $X+X' \in |-K_Y|$.  
Then $\Phi_{X,D}=0$.
\end{theo}
Note that this shows the log-Kobayashi conjecture for a very wide
class of log-$K$-trivial pairs.

The next proposition shows why this theorem is an easy consequence of
the existence on such pairs 
of log-$K$-autocorrespondences (\emph{i.e.} log-$K$-isocorrespondences
between a pair $(X,D)$ and itself, see definition \ref{logiso})
satisfying a certain dilation property.
\begin{prop}\label{auto-app}
Let $(X,D)$ be a pair composed by a smooth projective variety $X$, and
an effective divisor $D \subset X$, which is reduced and has normal
crossings, such that $K_X(D)$ is trivial. Let
$\eta_X$ be a generator of $\H^0(X,K_X(D))$. If there exists a
log-$K$-autocorrespondence $\Sigma$ of the pair $(X,D)$, such that for
a desingularization $\tau : \sigmat \rightarrow \Sigma$, and with the
notations (\ref{notations}) of definition \ref{Kcorresp}, one has
$$ f^* \eta_X = \lambda g^* \eta_X, $$
where $\lambda$ is a complex number with $|\lambda| \not= 1$, then
$\Phi_{X,D}=0$. 
\end{prop}

\proof
Let $\Omega_{X,D}$ be defined as
$$ \Omega_{X,D}=
(-1)^{\frac{n(n-1)}{2}} i^n \eta_X \wedge \bar{\eta}_X. $$
The dilation property satisfied by $\Sigma$ shows that
$$ f^* \Omega_{X,D} = |\lambda|^2 g^* \Omega_{X,D}, $$
while corollary \ref{isovolume} gives the equality of pseudo-volume
forms 
\begin{equation}\label{dil}
f^* \Phi_{X,D} = g^* \Phi_{X,D}.
\end{equation}
Now there exists a bounded, upper semi-continuous function $\chi$ on
$X$, such that $\Phi_{X,D} = \chi \Omega_{X,D}$.
$\chi$ has a maximum on $X$. Let $x$ be a point on $X$ where $\chi(x)$
is this maximum. Then take $\sigma \in \sigmat$ such that $f(\sigma)=x$,
and let $y=g(\sigma)$. Equality (\ref{dil}) eventually gives
$\chi(y)=|\lambda|^2 \chi(x)$. Since of course we can assume $|\lambda|
>1$ by symmetry, this shows that $\chi(x)=0$, and hence that $\chi=0$
as we wanted. 

\halmos

Now the proof of \ref{annulation} amounts to show the following
theorem. This is done in the remainder of this section.
\begin{theo}\label{existence}
If $(X,D)$ is a pair satisfying the condition of theorem
\ref{annulation}, then there exists a log-$K$-autocorrespondence
of $(X,D)$, such that (with the notations of proposition
\ref{auto-app}) 
\begin{equation}\label{fdiff}
m f^* \eta_X = -m' g^* \eta_X,
\end{equation}
where $m$ and $m'$ are distinct positive integers.
\end{theo}
In fact, by the following remark, it is enough to prove that there
exists a proper $0$-correspondence $\Sigma \subset X \times X$
(\emph{cf.} definition \ref{0corresp}) satisfying
relation (\ref{fdiff}). 
\begin{rem}\label{arg_Ktrivial}
Let $\Sigma \subset X \times X$ be a reduced closed analytic subset,
with both projections proper and generically finite. If $m f^* \eta_X
= -m' g^* \eta_X$ (with the notations of proposition
\ref{auto-app}), then $\Sigma$ is a log-$K$-autocorrespondence of
$(X,D)$.
\end{rem}
Indeed, $\eta_X$ is an everywhere non zero section of $K_X(D)$,
so the meromorphic $n$-form $f^* \eta_X$ (resp. $g^* \eta_X$) 
on $\sigmat$ has
divisor $R_f - f^*D$ (resp. $R_g - g^*D$), \emph{i.e.} it has zeroes
along $R_f$, and poles along $f^*D$. 
Equality (\ref{fdiff}) thus yields the equality of divisors  
$$ R_f - f^*D = R_g - g^*D $$
on the desingularization $\sigmat$, proving that $\Sigma$ is a
log-$K$-autocorrespondence of $(X,D)$.

\subsection{Geometric construction of self-correspondences on
  $K$-trivial pairs}

From now on, we let $Y$ be a smooth, rationally connected variety of
dimension $n+1$, and
$\bar{X} \subset Y$ be a reduced hypersurface in the anticanonical
linear class of $Y$ (\emph{i.e.} $\bar{X} \in |-K_Y|$), with smooth
locus $\sm(\bar{X}) \subset \bar{X}$. By adjunction, $\bar{X}$ has a
canonical sheaf $K_{\bar{X}}$, and it is trivial. We construct in this
section a $0$-correspondence $\Sigma \subset \bar{X} \times
\bar{X}$. 
In \cite{Kcorresp}, such  correspondences were constructed in the case
of smooth $\bar{X}$, and proved to be $K$-autocorrespondences.
The process is in some way analogous to the geometric
definition of addition on an elliptic curve in $\P^2$.
Of course, we have in mind the case when $\bar{X}=X \cup X'$ is the
reunion of two reduced hypersurfaces, $X$
smooth, and $D:=X \cap X'$ is reduced and has normal crossings, 
which will give theorem
\ref{existence} for the pair $(X,D)$, where $D$ is seen as a divisor
of $X$. 

Since $Y$ is rationally connected, there exists a rational curve $C_0
\subset Y$ satisfying the two following conditions. \\
(i) $C_0$ does not meet the singular part of $\bar{X}$, and $C_0 \cap
\bar{X} = mx_0 + m'y_0 + z_0$ as a divisor on $C_0$, where $x_0$ and
$y_0$ are distinct points of $\sm(\bar{X})$, $m$ and $m'$ are (fixed)
distinct positive integers, and $z_0$ is a reduced $0$-cycle on $C_0$,
disjoint from $x_0$ and $y_0$. \\
(ii) The deformations of the subscheme $C_0 \subset Y$ induce
arbitrary deformations of the $M$-jet of $C_0$ at two points of
intersection with $\bar{X}$, where $M:=\max(m,m')$. \\
This is given by the fact that $C_0$ can be chosen with arbitrarily
ample normal bundle, because $Y$ is rationally connected (see
\cite{KMM}). 

In addition, one chooses a hypersurface $W \subset \bar{X}$ containing the
$0$-cycle $z_0$. 
We denote by $|C_0|$ the space parametrizing all deformations of the
subscheme $C_0 \subset Y$.
Now define $\Sigma_{tot.} \subset \sm(\bar{X}) \times \sm(\bar{X})
\times |C_0|$ by
\begin{equation}\label{equations}\begin{array}{l}
\Sigma_{tot.} := \left\lbrace
(x,y,C)\ \ibox{s.t.}\ C\ \ibox{deformation of}\ C_0, \right.
\vspace{0.1cm}\hspace{3.3cm}
\\
\hspace{3.3cm} 
\left.
C\cap \bar{X} = mx+m'y+z,\ \mathrm{supp}(z) \subset W
\right\rbrace.
\end{array}\end{equation}
\begin{prop}\label{construction}
For a generic choice of $W$, the Zariski closure of $\Sigma_{tot.}$ in
$\bar{X} \times \bar{X} \times |C_0|$ has a unique $n$-dimensional
irreducible component through the point $(x_0,y_0,C_0)$. We call this
component $\Sigma'$. 
The Zariski closure of the projection of $\Sigma'$ in $\bar{X} \times
\bar{X}$ is irreducible of dimension $n$. Let us denote it
$\Sigma$. 
Then $\Sigma \subset \bar{X} \times \bar{X}$ is a $0$-correspondence.
\end{prop}
We call $f'$ and $g'$ the two morphisms $\Sigma' \rightarrow \bar{X}$
given by the projections of $\bar{X} \times \bar{X} \times |C_0|$ on
its first and second factors respectively. $\Sigma$ is an
$n$-dimensional component of the closure of $(f',g')(\Sigma')$. 

\vspace{0.2cm}\proof
The proposition essentially follows from the fact that, for a generic
choice of $W$, $\Sigma_{tot.}$ is smooth and of dimension $n=\dim
\bar{X}$ at the point $(x_0,y_0,C_0)$. This gives the existence of a
unique component $\Sigma'$ at once. On the other hand,
since the deformations of $C_0 \subset Y$ induce arbitrary
deformations of the $M$-jet of $C_0$ at $x_0$ and $y_0$, this implies
that the image $(f',g')(\Sigma') \subset \bar{X} \times \bar{X}$ has a
component of dimension $n$. 

So let us show that $\Sigma_{tot.}$ is smooth and of dimension $n$ near
$(x_0,y_0,C_0)$. We first study the Hilbert
scheme of $C_0 \subset Y$ at the infinitesimal neighbourhood of
$C_0$. Since $C_0$ is a rational curve, and its normal bundle
$N_{C_0/Y}$ is ample, we have $\h^1(C_0,N_{C_0/Y})=0$, and the Hilbert
scheme of $C_0 \subset Y$ is smooth, and of dimension
$$h^0(C_0,N_{C_0/Y}) = \chi(C_0,N_{C_0/Y}).$$
By the Riemann-Roch formula,
\begin{eqnarray*} 
\chi(C_0,N_{C_0/Y}) &=&
\deg(N_{C_0/Y}) + \mathrm{rg}(N_{C_0/Y}) (1-g) \\
&=& -K_Y \cdot C_0 + (2g-2) + n(1-g) = n-2 - K_Y \cdot C_0
\end{eqnarray*}
($g=0$ is the geometric genus of $C_0$).
Now to compute the dimension of $\Sigma_{tot.}$, we note that we impose to
the deformations $C$ of $C_0$ to meet $\bar{X}$ properly and in the
smooth locus $\sm(\bar{X})$ (this is open), and
to have intersection $C\cap \bar{X} = mx+m'y+z$, with $\mathrm{supp}(z)
\subset W$. This imposes at most $(m-1) + (m'-1) + \deg z$
conditions. In fact, for a generic choice of $W$, these conditions are
infinitesimally independent at the starting point $(x_0,y_0,C_0)$. We
deduce from this that $\Sigma_{tot.}$ is smooth at the neighbourhood of
$(x_0,y_0,C_0)$, and, since $\bar{X} \in |-K_Y|$, of dimension
$$
\left(-K_Y \cdot C_0 + n - 2\right) - \left(m+m'-2 +\deg z\right) 
= n -K_Y \cdot C_0 - \bar{X} \cdot C_0 = n.
$$

\halmos

\subsection{Realization as log-$K$-autocorrespondences} 

When $\bar{X} = X \cup X'$ is the union of two reduced hypersurfaces,
and if $X$ is smooth, then a section of $K_{\bar{X}}$ gives a
meromorphic $n$-form on the component $X$ by
restriction to the smooth locus of $\bar{X}$.
In this subsection, we show the following.
\begin{theo}\label{realisation}
Let $\Sigma \subset \bar{X} \times \bar{X}$ be as before (see proposition
\ref{construction}). We consider a desingularization $\tau : \sigmat
\rightarrow \Sigma$ and use the notations (\ref{notations}) of
definition \ref{Kcorresp}. Let
$\eta_{\bar{X}}$ be a generator of $\H^0(\bar{X},K_{\bar{X}})$ (recall
that $K_{\bar{X}}$ is 
trivial). We have the equality of meromorphic differential $n$-forms
on $\sigmat$ 
\begin{equation}\label{fdiff2}
mf^* \eta_{\bar{X}} + m' g^* \eta_{\bar{X}}=0.
\end{equation}
\end{theo}

Specializing this result to the case $\bar{X}=X \cup X'$ already
mentioned, we get the following corollary. Of course, it proves theorem
\ref{existence}, and because of proposition \ref{auto-app}, we get our
main theorem \ref{annulation} (\emph{cf.} subsection
\ref{struct_preuve}). 
\begin{coro}\label{corofin}
If $\bar{X}=X \cup X'$ is the union of two reduced hypersurfaces, $X$ smooth,
and $D:=X \cap X'$ reduced, then 
$\Sigma$ induces a $0$-correspondence $\Sigma_1 
:= \Sigma \cap (X \times X) \subset X \times X$.
If $\eta_X$ is a generator of $\H^0(X,K_{X}(D))$, then we have the
equality of meromorphic differential $n$-forms
$$mf_{1}^* \eta_X + m' g_{1}^* \eta_X=0$$
on a desingularization $\tilde{\Sigma}_{1} \rightarrow \Sigma_1$, where
$f_{1}$ and $g_{1}$ are the natural morphisms $\Sigma_1 \rightarrow
X$.
In particular, $\Sigma_1$ is a log-$K$-autocorrespondence of the pair
$(X,D)$. 
\end{coro}

\proof
In this case, $\Sigma \subset \bar{X} \times \bar{X}$ splits into four
components, contained in $X \times X$, $X\times X'$, $X' \times X$ and
$X' \times X'$ respectively.
Of course, the desingularization $\tau : \sigmat \rightarrow \Sigma$ 
splits as well. The first component
gives $\Sigma_1 \subset X \times X$, equipped
with a desingularization $\sigmat_1 \rightarrow \Sigma_1$ induced by
$\tau$. 

On the other hand, we have by adjunction
$$ \left.(K_{\bar{X}})\right|_{X} = \left.K_Y(\bar{X})\right|_{X} 
= \left.K_Y(X)\right|_X \otimes \O_{X}(X') = K_{X}(D). $$
In particular, since $K_{\bar{X}}$ is trivial, $K_X(D)$ is trivial as
well, and if $\eta_{\bar{X}}$ is a generator of 
$\H^0(\bar{X},K_{\bar{X}})$, then its restriction
$\eta_X := \left. \eta_{\bar{X}} \right|_{X}$ is a generator of 
$\H^0(X,K_{X}(D))$. Since $X$ is smooth, $\eta_X$ is a meromorphic
differential $n$-form, with polar divisor $D$.
It is now clear from theorem \ref{realisation} that we have the
equality of meromorphic $n$-forms
$$mf_{1}^* \eta_X + m' g_{1}^* \eta_X=0$$
on the desingularization $\tilde{\Sigma}_{1} \rightarrow
\Sigma_{1}$.
By remark \ref{arg_Ktrivial}, this implies that $\Sigma_1$ is a
log-$K$-autocorres\-pon\-dence of $(X,D)$.

\halmos

\vspace{0.2cm}
\noindent \textbf{Proof of theorem \ref{realisation}.}
We consider a desingularization $\tau': \tilde{\Sigma}' \rightarrow
\Sigma'$, with maps $\tilde{f} : \sigmat' \rightarrow \bar{X}$ and
$\tilde{g} : \sigmat' \rightarrow \bar{X}$. It is enough to show
(\ref{fdiff2}) on $\sigmat'$, that is
\begin{equation}\label{fdiff3}
m \tilde{f}^* \eta_{\bar{X}} = -m' \tilde{g}^* \eta_{\bar{X}}.
\end{equation}
The reason for this is that $(\Sigma,\mathrm{pr}_1,\mathrm{pr}_2)$ is the
Stein factorization of $(\Sigma',f',g')$ (and also of
$(\sigmat',\tilde{f},\tilde{g})$). 
In fact, since (\ref{fdiff3}) is an equality of meromorphic
differential forms, it 
is enough to prove it locally, and on a dense open subset of
$\sigmat'$. 

We construct three $0$-correspondences $\Gamma_{x,y}$, $\Gamma_z$ and
$\Gamma_C$ between $\Sigma'$ and $\bar{X}$. Each
of them is defined by its fiber over a generic point $\sigma =(x,y,C) \in
\Sigma'$. We describe these generic fibers as subschemes of $\bar{X}$. \\
-The fiber of $\Gamma_{x,y}$ over $\sigma$ is the $0$-cycle
$m.x + m'.y$. \\ 
-The fiber of $\Gamma_z$ over $\sigma$ is the $0$-cycle $C \cdot X -
m.x - m'.y$. For $\sigma$ generic, this is the part of the 
intersection $C \cap \bar{X}$ that lies on the hypersurface $W \subset
\bar{X}$. \\
-The fiber of $\Gamma_C$ over $\sigma$ is the $0$-cycle $C \cap
\bar{X}$. $\Gamma_C$ can also be defined globally as $\Crond
\cap (\Sigma' \times \bar{X}) \subset \Sigma' \times Y$, where $\Crond$ is
the universal curve over $\Sigma'$.
$$\xymatrix{
\Crond\ \ar@{^(->}[r] \ar[dr] & \Sigma' \times Y \ar[d] \\
& \Sigma'
}$$

Now, following an idea of Mumford (\cite{mumford}), we define the
pull-back of $\eta_{\bar{X}}$ by a $0$-correspondence $\Gamma \subset
\Sigma' \times \bar{X}$. 
In fact, since it is enough to prove (\ref{fdiff3}) locally and on a
dense open subset of $\sigmat'$, we only need to define the pull-back
of the restriction of $\eta_{\bar{X}}$ to a dense open
subset of $\bar{X}$. 
We thus define a smooth, dense open subset $\Sigma\zero \subset
\Sigma'$ in the following way (since it is smooth, $\Sigma\zero$ can
also be seen as a dense open subset of $\sigmat'$). 
We let $\pr_{12}$ be the projection $\bar{X} \times
\bar{X} \times |C_0| \rightarrow \bar{X} \times \bar{X}$. Then
$\Sigma\zero$ is the subset of $\Sigma' \cap \pr_{12}^{-1}(\sm(\bar{X})
\times \sm(\bar{X}))$, above which the projection $\Gamma_C \rightarrow
\Sigma'$ is finite. By restricting $\Sigma\zero$, one can assume that
it is smooth, and that it is a self-$0$-correspondence of a smooth,
dense open subset $X\zero \subset \bar{X}$.
We call $\Gamma\zero$ the restriction of $\Gamma$ to $\Sigma\zero
\times X\zero$, and $\varphi$ and $\psi$ the two projections of
$\Gamma\zero$ on $\Sigma\zero$ and $X\zero$ respectively. 
The key point here is the fact that $\varphi$ is proper.
Of course,
one can also assume $\Gamma\zero$ to be smooth.
$$\xymatrix@C=10pt{
& \Gamma\zero \ar[dl]_{\varphi} \ar[dr]^{\psi} & \\ 
\Sigma\zero && X\zero
}$$

We now define the pull-back by $\Gamma\zero$ of the restriction
$\eta\zero := \left.\eta_{\bar{X}}\right|_{X\zero}$. Since $X\zero
\subset \sm(\bar{X})$, $\eta\zero$ is a holomorphic differential
$n$-form on $X\zero$. Its pull-back $\psi^* \eta\zero$ is then a
holomorphic $n$-form on $\Gamma\zero$. Eventually, we let 
$$ (\Gamma\zero)^* \eta\zero := \varphi_* \psi^* \eta\zero, $$
where $\varphi_*$ is the trace map relative to $\varphi$. 
$(\Gamma\zero)^* \eta\zero$ is a holomorphic $n$-form on
$\Sigma\zero$. It extends to a meromorphic $n$-form on the smooth
variety $\sigmat'$~; we call it $\Gamma^* \eta_{\bar{X}}$.
The definition of the trace map $\varphi_*$ is classical. It goes as
follows. Let  
$\omega$ be a holomorphic $n$-form on $\Gamma\zero$. If $U \subset
\Sigma\zero$ is an open subset above which $\varphi$ is not ramified,
then $\varphi^{-1}U$ is the disjoint union of $d$ open subsets
$U_1,\ldots,U_d \subset \Gamma\zero$, which are all isomorphic to
$U$. On $U$, $\varphi_* \omega$ is defined as the sum
$$\left.\omega\right|_{U_1} + \cdots + \left.\omega\right|_{U_d}.$$
Then all these local definitions glue together to give a well defined
holomorphic $n$-form $\varphi_* \omega$ on $\Sigma\zero$, even if
$\varphi$ ramifies on $\Gamma\zero$. 

We now have all the definitions we need to end the proof. 
The equality (\ref{fdiff3}) is a simple consequence of the following
proposition.  
\begin{prop}\label{clef}
The holomorphic $n$-form $(\Gamma_{x,y}\zero)^* \eta\zero$ on
$\Sigma\zero$ vanishes. 
\end{prop}

To see why this implies (\ref{fdiff3}), we let $\Gamma_{f'}\zero \subset
\Sigma\zero \times X\zero$ and $\Gamma_{g'}\zero 
\subset \Sigma\zero \times X\zero$ denote the restrictions to
$\Sigma\zero \times X\zero$ of the
graphs of $f' : \Sigma' \rightarrow \bar{X}$ and $g' : \Sigma'
\rightarrow \bar{X}$ respectively. 
By construction, we have the
equality of $n$-cycles in $\Sigma\zero \times X\zero$
$$ \Gamma_{x,y}\zero = m \Gamma_{f'}\zero + m' \Gamma_{g'}\zero. $$
This implies that 
$(\Gamma_{x,y}\zero)^* \eta\zero = m f'^* \eta\zero + m' g'^*
\eta\zero.$ 
So proposition \ref{clef} yields the equality (\ref{fdiff3}) by
continuity. 

\halmos

\vspace{0.2cm}\noindent  \textbf{Proof of proposition \ref{clef}.} 
By construction, we have $\Gamma_{x,y} + \Gamma_z = \Gamma_C$ as
$n$-cycles in $\Sigma' \times \bar{X}$, 
and therefore $\Gamma_{x,y}\zero + \Gamma_z\zero =
\Gamma_C\zero$ in $\Sigma\zero \times X\zero$. 
The vanishing of proposition
\ref{clef} is a consequence of the two vanishings $(\Gamma_z\zero)^*
\eta\zero = (\Gamma_C\zero)^* \eta\zero = 0$ of holomorphic $n$-forms
on $\Sigma\zero$. 

The first one is given by the fact that $\Gamma_z$ 
is contained in $\Sigma' \times W \subset \Sigma' \times X$. Let $\nu
: \tilde{W} \rightarrow W\zero \subset X\zero$ be a desingularization
of $W\zero := W \cap X\zero$, and
consider the following $0$-correspondence obtained by base change.
$$\xymatrix@C=10pt{
& (\Gamma_z\zero)' \ar[dl]_{\tilde{\varphi}} \ar[dr]^{\psi}  \\
\Sigma\zero && \tilde{W} \ar[rr]^{\nu} && X\zero
}$$
One has $(\Gamma_z\zero)^* \eta\zero = (\Gamma_z\zero)'^* \nu^*
\eta\zero = 0$, because $\dim W < n$ and therefore $\nu^* \eta\zero =
0$.  

As for the second one, it is a consequence of proposition
\ref{factorisation} below. This proposition says that we have an
equality of meromorphic forms on $\sigmat'$
\begin{equation}\label{fact}
\Gamma_C^* \eta_{\bar{X}} = \Crond^* l_* \eta_{\bar{X}}.
\end{equation}
Here, $l_* : \H^0(\bar{X},K_{\bar{X}}) \rightarrow \H^1(Y,K_Y)$ is a
push-forward map induced by the inclusion $l: \bar{X} \subset Y$, and
$\Crond^* : \H^1(Y,K_Y) \rightarrow \H^0(\sigmat',K_{\sigmat'})$ is
given by the correspondence between smooth varieties
$\tilde{\Crond} \subset \sigmat' \times Y$, obtained by base change
from the universal family $\Crond \subset \Sigma' \times Y$. 
Now, since $Y$ is rationally connected, we have $\H^1(Y,K_Y)=0$. We
deduce from this that $l_* \eta_{\bar{X}} = 0$, and therefore that
$\Gamma_C^* \eta_{\bar{X}} = 0$. 
\emph{A fortiori}, we have
$(\Gamma_C\zero)^* \eta\zero = 0$ by restriction to $\Sigma\zero$.

\halmos

\vspace{0.2cm}
Let us now define precisely the various maps involved in formula 
(\ref{fact}). 
Let $l$ be the inclusion $\bar{X} \subset Y$. 
We define a push-forward morphism $l_*$, as the boundary morphism 
$\H^0(\bar{X},K_{\bar{X}}) \rightarrow \H^1(Y,K_Y)$ of the long exact
sequence associated to the Poincar{\'e} residue exact sequence
$$ 0 \rightarrow K_Y \rightarrow K_Y(\bar{X}) \rightarrow K_{\bar{X}}
\rightarrow 0.$$ 
On the other hand, we consider $\Crond \subset \Sigma' \times Y$,
which is the universal curve relative to the parametrizing variety
$\Sigma'$. We call $\Phi$ and $\Psi$ the two projections on $\Sigma'$
and $Y$ respectively. By base change, $\Crond$ gives a 
correspondence in $\sigmat' \times Y$. Desingularizing this new
correspondence, we obtain the following diagram, where
$\tilde{\Crond}$ is smooth, as well as $\sigmat'$ and $Y$. 
$$\xymatrix@C=10pt{
& \tilde{\Crond} \ar[dl]_{\tilde{\Phi}} \ar[dr]^{\tilde{\Psi}} & \\ 
\sigmat' && Y
}$$
This yields a Mumford pull-back $\Crond^* : \H^1(Y,K_Y) \rightarrow
\H^0(\sigmat',K_{\sigmat'})$, defined as the composition
$$\H^1(Y,K_Y) \xrightarrow{\tilde{\Psi}^*}
\H^1(\tilde{\Crond},K_{\tilde{\Crond}})
\xrightarrow{(\tilde{\Phi}^*)^T} \H^0(\sigmat',K_{\sigmat'}),$$
where the last morphism is obtained by Serre duality as the transpose
of the pull-back map $\Phi^* : \H^n(\sigmat',\O_{\sigmat'})
\rightarrow \H^n(\tilde{\Crond},\O_{\tilde{\Crond}})$.

\begin{prop}\label{factorisation}
We have 
$\Gamma_C^* \eta_{\bar{X}} = \Crond^* l_* \eta_{\bar{X}}$,
as an equality of meromorphic forms on $\sigmat'$.
\end{prop}

Again, this is an equality of meromorphic forms, so we can prove it
locally in the dense open subset $\Sigma\zero \subset \sigmat'$, where
it is an equality of holomorphic differential forms. 

Now the situation is the following. $\Crond$ is a family of curves
over $\Sigma'$. These curves are embedded in $Y$. $\Gamma_C$ is cut
out on $\Crond \subset \Sigma' \times Y$ by $\Sigma' \times \bar{X}$,
and since $\bar{X} \subset Y$ is a divisor, $\Gamma_C$ can be seen as
a family of divisors on the curves of the family $\Crond$. So basically,
proposition \ref{factorisation} is just a generalization of the
following result on curves, which is an application of the residue
theorem.  

\begin{lem}\label{courbes}
Let $C$ be a smooth curve, and $j : Z \subset C$ be an effective
divisor. We have 
a push-forward morphism $j_* : \H^0(Z,K_Z) \rightarrow \H^1(C,K_C)$,
taken from the long exact sequence associated to the Poincar{\'e} residue
exact sequence  
$$0 \rightarrow K_C \rightarrow K_C(Z) \rightarrow K_Z \rightarrow
0.$$ 
The following diagram is commutative.
$$\xymatrix{
\H^0(Z,K_Z) \ar[r]^{j_*} \ar[dr]_{\int_Z} & 
\H^1(C,K_C) \ar[d]^{\int_C} \\
& \C
}$$
\end{lem}

\proof
We first study the case when $Z$ is a
single point $x \in C$ with multiplicity $\mu$. It will then easily
extend to the general case, as will be seen from the proof. 
We begin by defining a $C^{\infty}$ function
$$
\rho : z\in \D \mapsto \exp\left(
  \frac{|z|^2}{|z|^2-1}\right) \in \R.
$$
It extends by $0$ on the complementary of the unit disk in $\C$ to a
$C^{\infty}$ function $\C \rightarrow \R$. It also satisfies
$\rho(0)=1$. 

We then consider a neighbourhood $U \subset C$ of $x \in C$, equipped
with a holomorphic coordinate $z$ centered at $x$, and such
that $U \cong \D$ \emph{via} $z$, and $Z \subset C$ is given by the
equation $z^{\mu}=0$. Let $\omega_Z \in \H^0(Z,K_Z)$. It writes 
\begin{equation}\label{meroinf}
\frac{a_{-\mu+1} + a_{-\mu+2}z + \cdots + a_{0}z^{\mu-1}}{z^{\mu-1}} 
\in \C[z]/(z^{\mu})\cdot \frac{1}{z^{\mu-1}},
\end{equation}
and we have $\int_Z \omega_Z = a_{0}$. 
The differential form
$$ \frac{\rho(z^{\mu})}{2\pi i} 
(a_{-\mu+1} + a_{-\mu+2}z + \cdots + a_{0}z^{\mu-1})
\frac{dz}{z^{\mu}} $$
is sent to $\omega_Z$ by the residue map. 
Its $\bar{\partial}$-differential 
\begin{equation}\label{diffdbarre}
\frac{1}{2\pi i} \frac{-\mu |z^{\mu-1}|^2}{(|z^{\mu}|^2-1)^2}
\exp\left( \frac{|z^{\mu}|^2}{|z^{\mu}|^2-1} \right)
\left( \frac{a_{-\mu+1}}{z^{\mu-1}} + \cdots + a_{0} \right)
d\bar{z} \wedge dz 
\end{equation}
thus represents $j_* \omega_Z$.
Note that in (\ref{diffdbarre}), the pole $z^{\mu-1}$ is eliminated by
the $|z^{\mu-1}|^2$ from the numerator. The differential form $j_*
\omega_Z$ is therefore a $C^{\infty}$ section of $\Omega^{1,1}_C$.
It is supported on $U \subset C$.
The integral of this form on $\D$ is
\begin{eqnarray*}
\frac{1}{2\pi i} \int_{\D}
\frac{-\mu (r^{\mu-1})^2}{(r^{2\mu}-1)^2} 
\exp \left( \frac{r^{2\mu}}{r^{2\mu}-1} \right)
\left( \frac{a_{-\mu+1} e^{-i(\mu-1)\theta}}{r^{\mu-1}} + \cdots + a_{0}
\right) 
(2ir dr \wedge d\theta) \\
= a_{0} \int_{0 \leqslant r \leqslant 1}
\frac{-2r^{\mu}}{(r^{2\mu}-1)^2} 
\exp \left( \frac{r^{2\mu}}{r^{2\mu}-1} \right)
(\mu r^{\mu-1}dr)
= a_{0}.
\end{eqnarray*}
This proves lemma \ref{courbes}.

\halmos

\vspace{0.2cm}
The next step of the proof of proposition \ref{factorisation} is to
extend lemma \ref{courbes} to the family of curves $\Crond$. In fact,
we have seen that it is enough to consider the smaller family
$\Crond\zero \rightarrow \Sigma\zero$.
$$\xymatrix@C=35pt{
\Gamma_C\zero \ar[r]^{(\mathrm{id}_{\Sigma\zero},l)} \ar[dr]_{\varphi} &
\Crond\zero \ar[d]^{\Phi} \\
& \Sigma\zero
}$$
We recall that $\varphi$ and $\psi$ are the natural morphisms from
$\Gamma_C$ to $\Sigma'$ and $X$ respectively, that $\Phi$ and $\Psi$
are the natural morphisms from $\Crond$ to $\Sigma'$ and $Y$, and that
$l$ is the inclusion $\bar{X} \subset Y$.
The following lemma is proved by using lemma \ref{courbes}, and the
fact that the push-forward map $\Phi_*$ is given locally above
$\Sigma\zero$ by the integration along the fibers of $\Phi$. 
\begin{lem}\label{famille}
The following diagram is commutative.
$$\xymatrix@C=35pt{
\H^0(\Gamma_C\zero,K_{\Gamma_C\zero})
\ar[r]^{(\mathrm{id}_{\Sigma\zero},l)_*} \ar[dr]_{\varphi_*} &
\H^1(\Crond\zero,K_{\Crond\zero}) \ar[d]^{\Phi_*} \\
& \H^0(\Sigma\zero,K_{\Sigma\zero})
}$$
\end{lem}
To be more precise, the vertical arrows of the preceding diagram are
given by $\varphi\zero_*$ and $\Phi\zero_*$ respectively, where
$\varphi\zero := \left. \varphi \right|_{\Gamma_C\zero}$ and
$\Phi\zero := \left. \Phi \right|_{\Crond\zero}$.

\proof
Again, we want to show an equality of holomorphic differential forms
on $\Sigma\zero$. This can be done locally. We thus choose an open
subset $V \subset \Sigma\zero$ over which the map $\Gamma_C\zero
\rightarrow \Sigma\zero$ is {\'e}tale (or rather the map
$(\Gamma_C\zero)_{red} \rightarrow \Sigma\zero$ is {\'e}tale). Then
$\varphi^{-1}(V)$ is a disjoint union of open subsets $V_j \subset
\Gamma_C\zero$, such that for every $j$ one has $(V_j)_{red} \cong
V$. 

One can assume that $V$ is equipped with a holomorphic system of
coordinates $v=(v_1, \ldots,$ $v_n)$, and that there exists a
neighbourhood $U_j$ of every $V_j \subset \Phi^{-1}(V)$ equipped with a
holomorphic system of coordinates $(v,t_j)$, such that $\Phi$ is simply
given locally by $(v,t_j) \mapsto v$, and that $V_j$ is given by the
local 
equation $t_j^{\mu_j}=0$ ($\mu_j$ is equal to $1$, $m$ or $m'$). The
coordinates $t_j$ give local parameters of the curves $C_{\sigma}$
($\sigma \in V$) at their points of intersection with $\bar{X}$.  

Let $\eta \in \H^0(\Gamma_C\zero,K_{\Gamma_C\zero})$. Its restriction
to $\varphi^{-1}(V)$ is a collection of meromorphic forms
$$ \eta_j = 
\left( a_{j,-\mu_j+1}(v) + a_{j,-\mu_j+2}(v) t_j + \cdots + a_{j,0}(v)
  t_j^{\mu_j-1} \right)
\frac{dv_1 \wedge \ldots \wedge dv_n}{t_j^{\mu_j-1}} $$ 
on the open subsets $V_j$, exactly as in (\ref{meroinf}).
The map $\varphi_*$ is locally defined as the trace map, so $\varphi_*
\eta$ is given in $V$ by
$$ \varphi_* \eta = \left( \sum\nolimits_j a_{j,0}(v) \right) dv_1
\wedge \ldots \wedge dv_n. $$ 
On the other hand, $(\id,l)_* \eta$ is given in $\Phi^{-1}(V)$ by a
$\bar{\partial}$-closed, $C^{\infty}$ differential form of type
$(n+1,1)$, which vanishes outside from the neighbourhood 
$\bigcup_j U_j$ of $\Gamma_C\zero \subset \Crond\zero$. In the
neighbourhood $U_j$ of $V_j \subset \Phi^{-1}(V)$, it is given by
$$ \bar{\partial} \left(
\frac{\rho(t_j^{\mu_j})}{2\pi i}
\left( a_{j,-\mu_j+1}(v) + a_{j,-\mu_j+2}(v) t_j + \cdots + a_{j,0}(v)
  t_j^{\mu_j-1} \right) 
\frac{dt_j}{t_j^{\mu_j}} \wedge dv_1 \wedge \ldots \wedge dv_n 
\right), $$
exactly as in the proof of lemma \ref{courbes}.
Eventually, $\Phi_*$ is given by the integration along the fibers of
$\Phi$. To compute this, we let $\zeta$ be a holomorphic $n$-vector
field over $V$ 
(\emph{i.e.} a holomorphic section of $\bigwedge^n T_{\Sigma\zero}$
over $V$), given in coordinates by
\begin{equation}\label{nvect}
h(v) \frac{\partial}{\partial v_1} \wedge \ldots \wedge
\frac{\partial}{\partial v_n}. 
\end{equation}
We lift it to a $C^{\infty}$ $n$-vector field $\tilde{\zeta}$ on
$\Phi^{-1}(V)$. Using a partition of unity, we can assume that
in every $U_j$, $\tilde{\zeta}$ is simply given by the expression
(\ref{nvect}). At a point $\sigma \in V$ of coordinates $v$,
the inner product
$\left(\Phi_* (\id,l)_* \eta\right) (\zeta)_{\sigma}$ is by definition
$$ \int_{C_{\sigma}} \left( (\id,l)_* \eta \right) (\tilde{\zeta}), $$
where $C_{\sigma}$ is the fiber of $\Phi$ over $\sigma$.  
We let $Z_{\sigma} := (\Gamma_C\zero \cap C_{\sigma})$, and denote by
$j_{\sigma}$ the inclusion $Z_{\sigma} \subset C_{\sigma}$.
The inner product 
$\eta(\tilde{\zeta})_{Z_{\sigma}} := \sum_{z \in
  \mathrm{Supp}(Z_{\sigma})} \eta(\tilde{\zeta})_{z}$  
sits naturally in $\H^0(Z_{\sigma},K_{Z_{\sigma}})$. 
Now the $(1,1)$-form $\left( (\id,l)_* \eta \right) (\tilde{\zeta})$
restricted to $C_{\sigma}$ is precisely 
$(j_{\sigma})_* (\eta(\tilde{\zeta})_{Z_{\sigma}} ),$
so by lemma \ref{courbes}, we have
$$ \left(\Phi_* (\id,l)_* \eta\right) (\zeta)_{\sigma} =
\int_{C_{\sigma}} (j_{\sigma})_* (\eta(\tilde{\zeta})_{Z_{\sigma}}) =
\int_{Z_{\sigma}} \eta(\tilde{\zeta})_{Z_{\sigma}}
= (\varphi_* \eta) (\zeta)_{\sigma}. $$
This shows the equality of holomorphic differential forms 
$\Phi_* (\id,l)_* \eta = \varphi_* \eta$ on $V$.

\halmos

\vspace{0.2cm}
To conclude the proof of proposition \ref{factorisation}, it only
remains to show the following commutativity result.
\begin{lem}\label{commut}
The following diagram is commutative.
$$\xymatrix@C=35pt{
\H^0(\bar{X},K_{\bar{X}}) \ar[r]^{l_*} \ar[d]_{\psi^*} &
\H^1(Y,K_Y) \ar[d]^{\Psi^*} \\
\H^0(\Gamma_C\zero,K_{\Gamma_C\zero})
\ar[r]^{(\mathrm{id}_{\Sigma\zero},l)_*} &
\H^1(\Crond\zero,K_{\Crond\zero})
}$$
\end{lem}
In fact, the vertical arrows of this diagram are rather $(\psi\zero)^*$
and $(\Psi\zero)^*$, where $\psi\zero : \Gamma_C\zero \rightarrow
\bar{X}$ and $\Psi\zero : \Crond\zero \rightarrow Y$ are the
restrictions of $\psi$ and $\Psi$ to $\Gamma_C\zero$ and $\Crond\zero$
respectively. 

\proof
It is clear from the definitions that $\psi\zero =
\left.\Psi\zero\right|_{\Gamma_C\zero}$. On the other hand, we have
$\Psi^* \bar{X} = \Gamma_C\zero$
as an equality of divisors on $\Crond\zero$. This shows that we have a
morphism of short exact sequences as follows.
$$\xymatrix{
0 \ar[r] & \Psi^* K_Y \ar[r] \ar[d]^{\Psi^*} &
\Psi^* K_Y(\bar{X}) \ar[r] \ar[d]^{\Psi^*} &
\psi^* K_{\bar{X}} \ar[r] \ar[d]^{\psi^*} & 0 \\
0 \ar[r] & K_{\Crond\zero} \ar[r] &
K_{\Crond\zero}(\Gamma_C\zero) \ar[r] &
K_{\Gamma_C\zero} \ar[r] & 0
}$$
It yields a morphism between the associated long exact sequences in
cohomology. In particular, the following diagram is commutative.
$$\xymatrix{
\cdots \ar[r] & \H^0(\Gamma_C\zero,\psi^* K_{\bar{X}}) \ar[r] \ar[d] & 
\H^1(\Crond\zero,\Psi^* K_Y) \ar[d] \ar[r] & \cdots \\
\cdots \ar[r] & \H^0(\Gamma_C\zero,K_{\Gamma_C\zero}) \ar[r] &
\H^1(\Crond\zero,K_{\Crond\zero}) \ar[r] & \cdots
}$$
Up to the shrinking of $\Sigma\zero \subset \Sigma$, we can assume
that both $\psi\zero$ and $\Psi\zero$ are smooth, so we have 
$\H^0(\Gamma_C\zero,\psi^* K_{\bar{X}}) \cong
\H^0(\bar{X},K_{\bar{X}})$ and $\H^1(\Crond\zero,\Psi^* K_Y) \cong
\H^1(Y,K_Y)$, and the lemma is proved.

\halmos

Eventually, lemmas \ref{famille} and \ref{commut} give the
commutativity of the left-hand side square and of the right-hand side
triangle in the following diagram.
$$\xymatrix{
\H^0(\bar{X},K_{\bar{X}}) \ar[r]^{\psi^*} \ar[d]_{l_*}
    \ar@/^2pc/[rr]^{(\Gamma_C\zero)^*} &
\H^0(\Gamma_C\zero,K_{\Gamma_C\zero}) \ar[r]^{\varphi_*}
    \ar[d]_{(\id_{\Sigma\zero},l)_*} & 
\H^0(\Sigma\zero,K_{\Sigma\zero}) \\
\H^1(Y,K_Y) \ar[r]_{\Psi^*} \ar@/_3.5pc/[urr]_{(\Crond\zero)^*} &
\H^1(\Crond\zero,K_{\Crond\zero}) \ar[ur]_{\Phi_*}
}$$
On the other hand, it follows from the definitions of
$(\Gamma_C\zero)^*$ and $(\Crond\zero)^*$ that the bended arrows
commute with the rest of the diagram. We thus have the equality of
holomorphic forms on $\Sigma\zero$
$$ (\Gamma_C\zero)^* \eta\zero = (\Crond\zero)^* l_* \eta_{\bar{X}},
$$  
which yields by continuity the desired equality of meromorphic
differential forms on $\sigmat'$. This ends the proof of proposition 
\ref{factorisation}. 
Theorem \ref{realisation} is thus completely proved.

\textsc{Thomas Dedieu, Universit{\'e} Paris 6 (Pierre et Marie Curie),
  Institut de Math{\'e}\-matiques de Jussieu, {\'E}quipe de Topologie et
  G{\'e}om{\'e}trie Alg{\'e}briques, 175 rue du Chevaleret, 75013 Paris, France}

\textit{E-mail address~:} \texttt{dedieu@math.jussieu.fr}

\end{document}